\documentstyle[12pt]{article}
\begin{document}
\newtheorem{Def}{Definition}[section]
\newtheorem{thm}{Theorem}[section]
\newtheorem{lem}{Lemma}[section]
\newtheorem{rem}{Remark}[section]
\newtheorem{prop}{Proposition}[section]
\newtheorem{cor}{Corollary}[section]
\title
{Compactness of solutions to the Yamabe problem. II
}
\author{ YanYan Li
 \\ Department of Mathematics\\ 
Rutgers University\\ 110 Frelinghuysen Road\\
Piscataway, NJ 08854\\ \\
Lei Zhang \\
Department of Mathematics\\
University of Florida\\
358 Little Hall\\
Gainesville, FL 32611-8105}
\date{}
\maketitle
\input { amssym.def}

\bigskip
\bigskip
\bigskip

\setcounter{section}{0}

\section{Introduction}

Let $(M^n, g)$ be a
compact, smooth, connected Riemannian manifold (without boundary)
of dimension $n\ge 3$.  The Yamabe conjecture has been proved
through the works of Yamabe \cite{Y}, Trudinger \cite{T}, 
Aubin \cite{A}
and Schoen \cite{S}:  There exist constant scalar curvature metrics
on $M$ which
are pointwise conformal to $g$. 

Consider the Yamabe equation 
and its sub-critical 
approximations:
\begin{equation}
-L_g u
=n(n-2)u^p,\qquad u>0,\quad \mbox{on}\quad M,
\label{Y0}
\end{equation}
where $1<p\le \frac{n+2}{n-2}$, 
$L_g = \Delta_g- c(n)R_g $, 
 $\Delta_g$ is the Laplace-Beltrami operator 
associated with $g$, $R_g$ is the
scalar curvature of $g$, and $c(n)=\frac { (n-2) }{ 4(n-1) }$.

Let 
$$
{\cal M}_p=\{ u\in C^2(M)\ |\ u\
\mbox{satisfies}\ (\ref{Y0})\}.
$$

Schoen initiated the investigation  of the
compactness of ${\cal M}_p$
and proved the following remarkable result in 1991, see \cite{S2},
under the assumption that  $(M,g)$ is locally conformally
flat and is not conformally diffeomorphic to
standard spheres:
For any $1<1+\epsilon\le p\le \frac {n+2}{n-2}$ and 
for any non-negative integer $k$,
\begin{equation}
\|u\|_{ C^k(M,g) } \le C,\qquad \forall\ u\in {\cal M}_p,
\label{compact}
\end{equation}
where $C$ is some constant depending only on
$(M,g)$, $\epsilon$ and $k$.  He also announced in the same
paper the same result for general manifolds, without the
locally conformally flat assumption.
The proof of this claim has not been made available.
For general manifolds of dimension $n=3$, a proof was given by Li and Zhu
in \cite{LiZhu}; while for $n=4$, a combination of
the results of Li and
Zhang \cite{LZ1} and Druet \cite{D1} yields a proof, with the $H^1$ bound
given in \cite{LZ1} and the $L^\infty$ bound under the assumption
of an $H^1$ bound
given in \cite{D1}.

Let $W_g$ denote the Weyl tensor (see for instance
\cite{A2} for the definition),
we consider the following two cases:

$1^\circ$.   $3\le n\le 7$,

$2^\circ$.   $n\ge 8$ and
$|W_g|+|\nabla W_g|>0$ on $M$.

\begin{thm}  Let $(M^n,g)$ be a 
compact, smooth, connected Riemannian manifold  which is not locally 
conformally flat.  
Assume either  $1^\circ$ or $2^\circ$.
Then, for any  $1<1+\epsilon\le p\le
\frac{n+2}{n-2}$ and for any non-negative integer $k$, (\ref{compact}) holds
for 
some constant $C$ depending only on
$(M^n,g)$, $\epsilon$ and $k$.
\label{thm3}
\end{thm}

\begin{rem}
The proof of Theorem \ref{thm3} in Case 1$^\circ$ makes
use of the deep positive mass theorem of Schoen and
Yau in \cite{SY0}.
If we assume the positive mass theorem for $n=8,9$, then our proof 
yields 
the conclusion of Theorem \ref{thm3} on any non-locally 
conformally flat
Riemannian manifolds of dimension $n=8,9$.
\label{rem0-1}
\end{rem} 

\begin{rem}   Theorem
\ref{thm3} in the case $n\le 4$ was, as mentioned earlier, already known.
   Theorem \ref{thm3} in the
case $n=5,6, 7$ as well as in the case $n\ge 8$ under
$|W_g|>0$ on $M$ was announced in November
2003 by the first author in his talk at the Joint Analysis
Seminar in Princeton University and also 
announced  in our note \cite{LZ2}.
These were also independently proved by Marques \cite{M}.
The case $n= 5$ was proved independently by Druet \cite{D2}.
 Theorem \ref{thm3} was announced by the first author 
 at the international 
conference in honor of Haim Brezis$'$s sixtyth birthday
 in Paris,  June 9-13, 2004. 
 Some further results for dimensions $n\ge 10$ will be given 
in a forthcoming paper \cite{LZ3}.
\label{0-2}
\end{rem}

Since 
the first eigenvalue of $-L_g$ is positive,  multiplying (\ref{Y0}) by
$u$ and integrating by parts on $M$ lead to  
 $\max_{M}u\ge C^{-1}$  for some positive constant $C$ depending only  
$(M,g)$ and $\epsilon$.  
On the other hand, once we know that
$u\le C$ on $M$ for some  $C$ depending only on
$(M,g)$ and $\epsilon$, an application of the
Harnack inequality yields 
$u\ge C^{-1}\max_{M}u$ on $M$
for some  $C$ depending only on
$(M,g)$ and $\epsilon$.
A consequence of Theorem \ref{thm3} is, by some arguments in \cite{S2}
and \cite{Li1},
the following
\begin{cor}
Under the hypothesis of Theorem \ref{thm3}, there exists some
constant $\bar C$, depending only on $(M,g)$ and
$\epsilon$, such that for all $1<1+\epsilon\le p\le \frac{n+2}{n-2}$
and for all $C>\bar C$,
$$
\mbox{deg} \left( v-n(n-2) (L_g)^{-1}(v^p), {\cal O}_C, 0\right)
=-1,
$$
where
${\cal O}_C:=\{v\in C^{2,\alpha}(M)\ :\
1/C<v<C, \|v\|_{  C^{2,\alpha}(M) }<C\}$,
$0<\alpha<1$, and $deg$ denotes the Leray-Schauder degree.
\label{cor0-1}
\end{cor}

For the Leray-Schauder degree theory, see for instance
\cite{N}.

Much of this paper  is devoted to the fine analysis of blow up
 solutions of
(\ref{Y0}).  
In fact we establish the following local version of such estimates.
Let 
 $\Omega\subset M$ be
an open connected subset of $M$, and let
$\Omega_\epsilon:=\{P\in \Omega\ |\ dist_g(P, \partial\Omega)>\epsilon)$.
For  $Q\in \Omega$ and $\mu>0$, let
$$\xi_{Q,\mu}(P)=
(\frac{\mu}{1+\mu^2\mbox{dist}_g(P,Q)^2})^{\frac{n-2}{2}},
\quad P\in \Omega. $$

We are interested in solutions  of
\begin{equation}
-L_gu=n(n-2)u^p,\qquad u>0,\quad \mbox{in}\ \Omega.
\label{YY3}
\end{equation}

\begin{thm}
Let $(M^n,g)$ be a compact, smooth, connected                                     Riemannian manifold of dimension $n\ge 10$,
and let  $\Omega\subset M$ be
an open connected subset, $1<1+\epsilon\le p\le \frac{n+2}{n-2}$.
Suppose that
 $u$ is a smooth  solution of
 (\ref{YY3}) satisfying,
 for some $\bar P\in \Omega$ and some constant
 $\bar b\ge 1$, that
 \begin{equation}
 \nabla u(\bar P)=0,
 \qquad 1\le  \sup_\Omega u\le \bar b u(\bar P).
 \label{Y1}
 \end{equation}
 Then, for any $\epsilon'>0$, 
 \begin{equation}
 |W_g(\bar P)|_g+|\nabla _g W_g (\bar P)|_gu(\bar P)^{-\frac{2}{n-2}}
 \le C
 u(\bar P)^{ -\frac 4{n-2}+\epsilon'}.
 \label{mm2}
 \end{equation}
where $C$ is
  some positive constant 
  depending only 
   on                                
    $\epsilon$, $\epsilon'$,
    $ dist_g(\bar P, \partial \Omega)$, $\bar b$,
     a positive lower bound for the injectivity radius
      of $(M,g)$,  a positive lower bound of the first
       eigenvalue of $-L_g$ on $\Omega$ with zero
        Dirichlet boundary condition, and a positive upper bound
	of the norm of the curvature tensor of $(M,g)$ together with
	its covariant derivatives up to the eighth order.
\label{thm10}
\end{thm}
\begin{rem}  Theorem \ref{thm3} in the case 
$n\ge 10$ follows immediately from Theorem \ref{thm10}.
\label{rem10}
\end{rem}

The estimates we establish for dimensions $n\le 9$ in the next theorem is
much stronger than that for $n\ge 10$.  This is the reason that
 Theorem  \ref{thm3} for $n\le 7$,
as well as for dimensions $n=8,9$ as mentioned in 
Remark \ref{rem0-1}, hold without any assumption on the Weyl tensor.
\begin{thm}
Let $(M^n,g)$ be a compact, smooth, connected
Riemannian manifold of dimension $3\le n\le 9$,
and let  $\Omega\subset M$ be
an open connected subset, $1<1+\epsilon\le p\le \frac{n+2}{n-2}$.
Suppose that
 $u$ is a smooth  solution of
(\ref{YY3}) satisfying (\ref{Y1})
for some $\bar P\in \Omega$ and some constant
$\bar b\ge 1$.
Then, for any $\delta>0$,
there exist some positive constant
$C$ and some positive integer $m$, which
 depend only
 on 
 $\epsilon$, $ dist_g(\bar P, \partial \Omega)$, $\bar b$, 
$\delta$,   
 a positive lower bound for the injectivity radius
 of $(M,g)$, a positive lower bound for  the first
 eigenvalue of $-L_g$ on $\Omega$ with zero
 Dirichlet boundary condition, and a positive upper bound
of the norm of the curvature tensor of $(M,g)$ together with
its covariant derivatives up to the eighth order, 
 and there exist local maximum points
${\cal S}:=\{P_1, \cdots, P_m\}\subset \Omega_{2\delta}$ of $u$,
such that 
\begin{equation}
\mbox{dist}_g(P_i,P_j)\ge \frac 1{C},\qquad \forall i\neq j,
\label{ee7}
\end{equation}
\begin{equation}
\frac 1Cu(P_i)\le u(P_j)\le Cu(P_i),\qquad \forall i,j,
\label{2-5}
\end{equation}
\begin{eqnarray}
|W_g(P_i)|_g
&\le
\left\{
\begin{array}{rl}
 \frac C{\sqrt{\log u(P_i)}},& \mbox{if}\ n=6,\\
 C u(P_i)^{ -\frac {n-6}{n-2} },& \mbox{if}\ 7\le n\le 9, 
 \end{array}
 \right.
 &\qquad
 \forall\ i
 \label{nn1new}
 \end{eqnarray}
 \begin{eqnarray}
 |\nabla W_g(P_i)|_g&\le \left\{
 \begin{array}{ll} \frac C {\sqrt{\log u(P_i)}},& \mbox{if}\ n=8,\\
 C u(P_i)^{ -\frac {n-8}{n-2} },& \mbox{if}\ n= 9.
 \end{array}
 \right.&
 \qquad\forall\ i
 \label{nn2new}
 \end{eqnarray}
 \begin{equation}
\frac 1C\sum_{l=1}^{m}\xi_{P_l,u(P_l)^{ \frac 2{n-2} }}(P)\le u(P)\le
C\sum_{l=1}^{m}\xi_{P_l,u(P_l)^{ \frac 2{n-2} }}(P),
\quad \forall P\in \Omega_{4\delta}.
\label{2-6}
\end{equation}
and, for $|\alpha|=0,1,2$ and $P\in \Omega_{4\delta}$,
\begin{eqnarray}
&&|\nabla^{\alpha}_g(u-\sum_{l=1}^m\xi_{P_l,u(P_l)^{\frac{2}{n-2}}})(P)|
\nonumber\\
&\le& 
\left\{
\begin{array}{ll}
Cu(\bar P)^{-1+\frac{2|\alpha |}{n-2}}
(1+u(\bar P)^{\frac{2}{n-2}}
\mbox{dist}(P,{\cal S}))^{-|\alpha |},& \mbox{if}\ n=3,4, 5,\\
 C(\epsilon')u(\bar P)^{-1+\frac{2|\alpha |+2\epsilon'}{n-2}}
(1+u(\bar P)^{\frac{2}{n-2}}
\mbox{dist}(P,{\cal S}))^{-\epsilon' -|\alpha |}, & \mbox{if}\ n=6,\\
 Cu(\bar P)^{1+\frac{2|\alpha |-8}{n-2}}
(1+u(\bar P)^{\frac{2}{n-2}}
\mbox{dist}(P,{\cal S}))^{6-n-|\alpha |},  & \mbox{if}\ n=7,8,9.
\end{array}
\right.
\label{abc}
\end{eqnarray}
\label{thm0}
\end{thm}

Estimates (\ref{ee7}), (\ref{2-5}) and (\ref{2-6}) on locally
conformally flat manifolds, without the assumption (\ref{Y1}),
 were established  
by Schoen in \cite{S1}.
Analogues of such local results do not hold
for many other problems, including
Harmonic maps and  Ginzburg-Landau vortices,
with similar loss of compactness
 --- similar in the sense that
 the energy is quantized when a sequence 
 of solutions blows up.
It is interesting to note that
the  analogue of  
 (\ref{ee7}) 
 in such a local setting
(under the assumption (\ref{Y1}))
fails even for the equation
$-\Delta u=V e^u$ in 
 dimension $n=2$,
 as demonstrated by
Chen in \cite{ChenX}.  On the other hand, it was proved 
by the first author in
 \cite{Li3} that analogues of  (\ref{ee7}), (\ref{2-5}) and (\ref{2-6}) 
hold for solutions of  such equations  
defined globally 
on compact Riemannian surfaces.

The proof of the Yamabe conjecture
(\cite{Y}), \cite{T}, \cite{A} and \cite{S}),
 a milestone in the studies of nonlinear
elliptic equations, concerns the existence of a minimizer 
of some functional with lack of compactness.
Many further studies have been devoted to 
related  critical exponent equations 
which address important issues
including those concerning non-minimal solutions
or approximate solutions to such equations, 
see for instance Brezis and Nirenberg (\cite{BN}) and Bahri and Coron
(\cite{BC1} and \cite{BC2}).  These studies have led to 
different proofs of the Yamabe conjecture
in the case $n\le 5$ and in the case
$(M,g)$ is locally conformally flat, see
Bahri and Brezis \cite{BB} and Bahri \cite{B}.
For
 the
 Nirenberg problem and the Yamabe problem on
 manifolds with boundary,
 which are related to the Yamabe problem,
 compactness 
of solutions
has been studied 
by Schoen (\cite{S1}), 
Schoen and Zhang (\cite{SZ}), Chang, Gursky and Yang
(\cite{CGY}),  Li (\cite{Li1}) and \cite{Li2}),
Han and Li (\cite{HL}),
 Chen and Lin (\cite{ChenLin2}), 
Felli and Ould Ahmedou (\cite{FO}), and Escobar and Garcia
\cite{EG}.
  Much of the analysis in these works
can be made purely local.     
One way to achieve this
is a Harnack type inequality
of Schoen (\cite{S1}): For $n\ge 3$,
let $u$ be
a smooth positive solution of
$$
-\Delta u=u^{ \frac{n+2}{n-2} },
\qquad\mbox{in}\ B_4\subset \Bbb R^n,
$$
then
\begin{equation}
\sup_{B_1}u
\cdot \inf_{B_2}u\le
C(n).
\label{har}
\end{equation}
A consequence of this is, as proved in \cite{S1},
$\int_{ B_1} u^{ \frac{2n}{n-2} }\le C(n).
$
Analogous results, extensions, as well as
different proofs of (\ref{har}) can be found in
\cite{Siu},   \cite{Tian},
 \cite{S1},  \cite{BM}, 
 \cite{BLS},
 \cite{LS},
  \cite{ChenLin1},
  \cite{ChenLin0}, 
\cite{Li3},
\cite{BT0}, \cite{BT1},
\cite{ChenLin3}, 
\cite{LZ0}, 
\cite{LL1},
 \cite{Ba}, \cite{PR},
\cite{ChenLin4}, 
\cite{LL2}, 
\cite{BCLT}, \cite{Ta},
\cite{Ta2}, and others. 
In particular, it was proved in our 
paper \cite{LZ1} that such Harnack type inequality
holds 
 on non-locally conformally flat Riemannian manifolds
of dimension $n=3,4$. Therefore
the conclusion of Theorem \ref{thm0}
holds in dimension $n=3,4$ without the assumption           
 (\ref{Y1}).

The main step in the proof of Theorem \ref{thm10}
and Theorem \ref{thm0} is 
to establish Theorem \ref{thm1}.  The proof of Theorem \ref{thm1}
is based on the method of moving planes, using the  ansatz of Schoen
in his proof of (\ref{har}) in \cite{S1} (see also  \cite{BLS},
\cite{ChenLin1}, \cite{ChenLin2},
  \cite{Li3} and \cite{LZ1} where such  ansatz was used).
  The method of moving planes has become a powerful tool
  in the study of nonlinear elliptic equations, see Alexandrov \cite{Al},
  Serrin \cite{Serrin}, 
  Gidas, Ni and Nirenberg \cite{GNN},   Berestycki and
  Nirenberg \cite{BeN}, and others.  
The main task in our proof of
Theorem \ref{thm1}
 is to construct  suitable auxiliary functions so
 that the method of moving planes can be applied.
 To do this we make use of numerous  results
 and methods from previous works, some of which are described below. 
We need to make use of the Liouville type
 theorem of Caffarelli, Gidas and Spruck in \cite{CGS}
 which identifies the limit of the rescaled blow-up
  sequence of solutions, and we need to establish
 strong enough convergence rate of the difference of the rescaled blow-up 
 sequence of solutions and its limit.
 Certain rate of convergence of the difference of the rescaled blow-up
  sequence of solutions and its limit
  was established by Chen and Lin
 in \cite{ChenLin2} for the  scalar curvature equations
in the Euclidean space, and we
 adapt their methods to
 establish iterated estimates on such convergence rate  
 on larger and larger balls with improved estimates
 after each iteration.
 For dimensions $n\ge 8$, the iterated estimates also yield
 stronger and stronger decay estimates on the Weyl
 tensor and its first covariant derivatives
 at points where the sequence of
 solutions blow up.  Here we also make use
 of a Pohozaev type identity as well as some
 properties of the conformal normal coordinates
 as established by Lee and Parker \cite{LP},
 Cao \cite{Cao},   G\"unther \cite{G}, 
 and Hebey and  Vaugon
 \cite{HV}.
 To construct auxiliary functions we also
 adapt the way of using the
 spherical harmonics
  by
    Caffarelli,
    Hardt and Simon
     in \cite{CHS}.
     The proof of Theorem \ref{thm3} is based on
     Theorem  \ref{thm10}-\ref{thm0}
     and the positive mass theorem of Schoen and Yau in 
     \cite{SY0}.
     Theorem  \ref{thm10}-\ref{thm0} provide
     strong enough pointwise  estimates
     for blow-up solutions as well as, for higher dimensions,
     strong enough decay estimates
      for the Weyl tensor and its first covariant 
      derivatives at the blow up points.
     These estimates allow us
     to use the positive mass theorem through the Pohozaev
     type identity as in Schoen  \cite{S0}.
     We note that the Pohozaev identity has been used
     by Arkinson and Peletier \cite{AP} and
     Brezis and Peletier \cite{BP} to obtain pointwise 
     estimates for blow up solutions of related critical exponent
     equations,  see also  \cite{S1}, \cite{Han2},  \cite{SZ},
     \cite{Li1}, \cite{Li2}, and others, for the extensive use of
     the Pohozaev  identity in establishing pointwise 
      estimates to blow up solutions of critical
       exponent equations.

Most of the works mentioned above concern the
compactness of solutions or
the fine pointwise analysis of blow up solutions to the Yamabe
equation and some related ones, which often 
yield the existence of solutions through the use
of degree theories.  There
have been many works on the existence of 
solutions to the Yamabe problem, the Nirenberg problem, 
and the Yamabe problem on manifolds with boundary, see
for instance  \cite{Mo}, 
 \cite{KW}, 
 \cite{C},
 \cite{DN},  \cite{Hong},
 \cite{ES},  \cite{CY1},
 \cite{CY2}, 
 \cite{CD},
 \cite{Chen},  \cite{H},
 \cite{He},
 \cite{CY3},
 \cite{E1},  \cite{E2},
 \cite{BE},
 \cite{CL}, 
 \cite{E3}, 
 \cite{AB},    
 \cite{AB2},
 \cite{ChL},
\cite{AGP}, 
 \cite{HL2} and  \cite{ALM}.
There have also been works on parabolic flows associated with
the Yamabe problem and the Yamabe problem on manifolds with boundary,
see for instance 
  \cite{Ye} , 
  \cite{SS} and
 \cite{Br}.

\bigskip

\noindent{\bf Acknowledgment.}\
{\it We thank H. Brezis, L. Nirenberg and S. Taliaferro for
encouragement and stimulating 
discussions.
Part of this paper was completed while the first author was
a visiting member at the Institute for Advanced Study in
Fall 2003.  He thanks J. Bourgain and IAS for providing him the
excellent environment, as well as
for providing him the financial support through
NSF-DMS-0111298.  Part of the work of the first author
is also supported by NSF-DMS-0100819 and
NSF-DMS-0401118.
}

\section{Main estimates}

Let $B_1\subset \Bbb R^n$, $n\ge 3$,  be the
unit ball centered at the origin, and let
$(a_{ij}(x))$ be a smooth,  $n\times n$ symmetric
positive definite matrix function, defined
on $B_1$,  satisfying
\begin{equation}
\frac 12 |\xi|^2\le a_{ij}(x)\xi^i\xi^j\le 2  |\xi|^2,
\qquad \forall\ x\in B_1, \ \xi\in \Bbb R^n,
\label{cond1}
\end{equation}
and, for some $\bar a>0$,
\begin{equation}
\|a_{ij}\|_{ C^{8}(B_1) }\le \bar a.
\label{cond2}
\end{equation}
Consider the Riemannian metric
\begin{equation}
g=a_{ij}(x)dx^idx^j
\label{2-1}
\end{equation}
on $B_1$, and consider 
\begin{equation}
-L_g u=n(n-2)u^p,\qquad u>0,\quad \mbox{on}\quad B_1.
\label{eq3}
\end{equation}

\begin{thm}
 Let $(B_1,g)$ be as  above 
 and let
 $u$ be  a solution of (\ref{eq3}), with $1<1+ \epsilon
\le p\le \frac {n+2}{n-2}$, satisfying, for some $\bar b\ge 1$,
\begin{equation}
\nabla u(0)=0, \qquad 1\le \sup_{ B_1}u\le \bar b u(0).
\label{normal}
\end{equation}
 Then  there exist some
positive constants $\delta$ and $C_0$, depending only on $n$,  
$\bar b$, $\epsilon$ and $\bar a$,  such
that
\begin{equation}
u(0)u(x)|x|^{n-2}\le
C_0, \qquad \forall\ 0<|x|\le\delta, \qquad \mbox{if}\ \ \ \  3\le n\le 9,
\label{eq4}
\end{equation}
\begin{equation}
|W_g(0)|_g\le
\left\{
\begin{array}{rl}
 \frac {C_0}{\sqrt{ \log u(0) }  },& \mbox{if}\ n=6,\\
C_0  u(0)^{ -\frac {n-6}{n-2} },& \mbox{if}\ 7\le n\le 9,
\end{array}
\right.
\label{nn1}
\end{equation}
\begin{equation}
|\nabla_g W_g(0)|_g\le
\left\{
\begin{array}{rl}
 \frac {C_0}{\sqrt{\log u(0)}},& \mbox{if}\ n=8,\\
C_0u(0)^{ -\frac {n-8}{n-2} },& \mbox{if}\ n= 9.
\end{array}
\right.
\label{nn2}
\end{equation}
If $n\ge 10$, then for all $\epsilon_1>0$, there exists $C(\epsilon_1)>0$ 
such that 
\begin{equation}
|W_g(0)|_g+|\nabla_g W_g(0)|_gu(0)^{ -\frac 2{ (n-2) } }
\le C(\epsilon_1) u(0)^{ -\frac 4{ (n-2) }+\epsilon_1 }.
\label{mm1}
\end{equation}
\label{thm1}
\end{thm}

\begin{rem} Theorem \ref{thm10} follows from 
(\ref{mm1}).
\end{rem}

We first prove
 Theorem \ref{thm1} for $p=\frac{n+2}{n-2}$.
We point out the changes needed for $p<\frac{n+2}{n-2}$
in Section 5.
  Suppose that the conclusion of  Theorem \ref{thm1} for $p=\frac{n+2}{n-2}$
  does not hold,
then for some $\bar a>0$, $\bar b\ge 1$, there exist
a sequence of Riemannian metrics $\{\tilde g_k\}$
of the form  (\ref{2-1}) that satisfy (\ref{cond1})
and (\ref{cond2}), and some solutions $u_k$ of
(\ref{eq3}), with $p=\frac {n+2}{n-2}$ and
with $g$ replaced by $\tilde g_k$,
satisfying (\ref{normal}), such that
one of the following happens:
\begin{equation}
\max_{ |x|<\frac 1k} \bigg( 
u_k(0) u_k(x) |x|^{n-2} \bigg)\ge k,
\label{5-0}
\end{equation}
\begin{equation}
|W_{\tilde g_k}(0)|_{\tilde g_k}>
\left\{
\begin{array}{rl}
 \frac {k}{\sqrt{\log u_k(0)}},& \mbox{if}\ n=6,\\
k u_k(0)^{ -\frac {n-6}{n-2} },& \mbox{if}\ 7\le n\le 9,
\end{array}
\right.
\label{555}
\end{equation}
\begin{equation}
|\nabla_{\tilde g_k}W_{\tilde g_k}(0)|_{\tilde g_k}>
\left\{
\begin{array}{rl}
 \frac {k}{\sqrt{\log u_k(0)}},& \mbox{if}\ n=8,\\
ku_k(0)^{ -\frac {n-8}{n-2} },& \mbox{if}\ n= 9,
\end{array}
\right.
\label{nn9}
\end{equation}
or, for some  $\epsilon_5>0$ independent of $k$,  
\begin{equation}
|W_{\tilde g_k}(0)|_{ \tilde g_k}
+|\nabla _{\tilde g_k} W_{\tilde g_k}(0)|_{ \tilde g_k}u_k(0)^{-\frac{2}{n-2}}
>k u_k(0)^{ -\frac 4{(n-2) }+\epsilon_5 },
\qquad \mbox{if}\ n\ge 10.
\label{mm4}
\end{equation}
We will simply use $g$ to denote
$\tilde g_k$.

Let  $\bar P$ be a point on
$(M,g)$, it was proved in  \cite{LP},
together with some improvement in \cite{Cao} and \cite{G},  that
there exists some function  $\varphi$ (with control)
near $\bar P$  such that the conformal metric
$\tilde g=e^\varphi g$ satisfies, in $\tilde g-$normal coordinates
$\{x^1, \cdots, x^n\}$ centered at $\bar P$
$$
\det(\tilde g_{ij})=1.
$$
Such coordinates are called conformal normal
coordinates. As well known, we may  assume that we work
in conformal normal coordinates.
In conformal normal coordinates (we write $g_{ij}$ instead of $\tilde g_{ij}$),
we have, at $x=0$,
\begin{equation}
R_{ij}=0, \  R_{, i}=0,
\  Sym_{ijk}R_{ij,k}=0,\   \Delta_g R= -\frac 16|W|^2,
\label{conformal1}
\end{equation}
where $R_{ijkl}$ denotes the curvature tensor evaluated at $0$,
$R_{ij}$ denotes the Ricci curvature tensor at $0$,
$R_{ijkl, p}$ denotes covariant derivative of
the curvature tensor at $0$, etc.,
repeated indices denote summation over the indices,  and
$$
Sym_{p_1\cdots p_m}A_{p_1\cdots p_m}:=\sum_{\sigma}
A_{p_{\sigma(1)} \cdots p_{\sigma(m)} }
$$
where $\sigma $ runs through all permutations of $1, 2, \cdots, m$.

We make a conformal change of the metric
$\hat g=e^\varphi g_k$
and let  $\{z^1, \cdots, z^n\}$  be the conformal  normal coordinates
centered at the origin.  
After the conformal change, $u_k$ becomes $\hat u_k=e^{ \frac {n-2}4
\varphi} u_k$.  As well known all relevant properties
of $u_k$ hold for $\hat u_k$, and we simply
assume that $g_{ij}(z)dz^idz^j$ is already in conformal
normal coordinates. In local coordinates,
\begin{eqnarray}
g_{pq}(x)&=&\delta_{pq}+\frac 13R_{pijq}x^ix^j+\frac 16R_{pijq,k}x^ix^jx^k
\nonumber \\
&&+(\frac 1{20}R_{pijq,kl}+\frac 2{45}R_{pijm}R_{qklm})x^ix^jx^kx^l+O(r^5).
\nonumber
\end{eqnarray}
In conformal normal coordinates, write
$$
\Delta_g=\frac 1{\sqrt{g}}\partial_i(\sqrt{g}g^{ij}\partial_j)
=\Delta +b_i\partial_i +d_{ij}\partial_{ij},
$$
where $(g^{ij})$ denotes the inverse matrix of $(g_{ij})$,
$\partial_i=\frac {\partial }{\partial z^i}$,
 $\partial_{ij}=\frac {\partial ^2}{ \partial z^i \partial z^j}$,
$\Delta=\sum_{i=1}^n \frac {\partial ^2}{ \partial z^i \partial z^i}$,
\begin{eqnarray}
&&b_i(x)= \partial_j g^{ij}(x)\nonumber\\
&=&-\frac 16 R_{ia,b}x^ax^b-\frac 16R_{iabp,p}x^ax^b-(\frac 1{20}
R_{ia,bc}-\frac 1{15}R_{ipad}R_{pbcd}
\nonumber\\
&&-\frac 1{15}R_{iapd}R_{pbcd}+\frac 1{20}R_{iabp,pc}
+\frac 1{20}R_{iabp,cp})x^ax^bx^c+O(r^4),
\nonumber
\end{eqnarray}
and
\begin{eqnarray}
d_{ij}(x)=g^{ij}-\delta_{ij}=-\frac 13R_{ipqj}x^px^q-\frac 16R_{ipqj,k}x^px^qx^k
-(\frac 1{20}R_{ipqj,kl}\nonumber \\
-\frac 1{15}R_{ipqm}R_{jklm})x^px^qx^kx^l+O(r^5).\nonumber
\end{eqnarray}

By (\ref{normal}) and (\ref{5-0}),
 $M_k:=u_k(0)\to \infty$.   
Write $(g_k)_{ij}(y)=g_{ij}(M_k^{ -\frac 2{n-2} }y)dy^idy^j$,
then
$$\Delta_{g_k}=\Delta+\bar b_i\partial_i+\bar d_{ij}\partial_{ij},
$$
where
$$
 \bar b_i(y)=M_k^{-\frac{2}{n-2}}b_i(M_k^{-\frac{2}{n-2}}y),
\quad \bar d_{ij}(y)=d_{ij}(M_k^{-\frac{2}{n-2}}y).
$$

Let
\begin{equation}
v_k(y):=M_k^{-1}u_k(M_k^{-\frac{2}{n-2}}y),
\label{vk}
\end{equation}
$$
c(x)=c(n)R_g(x), \quad \mbox{and} \quad
 \bar c(y)=c(n)R_g(M_k^{-\frac 2{n-2}}y)M_k^{-\frac 4{n-2}}.
$$
Then
\begin{equation}
|\bar b_i(y)|=O(1) M_k^{-\frac{6}{n-2}} |y|^2,
 \
|\bar d_{ij}(y)|=O(1)M_k^{-\frac{4}{n-2}}|y|^2,
 \
\bar c(y)=O(1)M_k^{ -\frac 8{n-2} }|y|^2.
\label{rough}
\end{equation}

The rescaled function $v_k$
 satisfies
\begin{equation}
\left\{
\begin{array}{ll}
\Delta_{g_k}v_k(y)-\bar cv_k(y)+n(n-2)v_k(y)^{\frac{n+2}{n-2}}=0,
\quad |y|\le \frac 12M_k^{\frac{2}{n-2}},
  \\
1=v_k(0)\ge (\bar b^{-1}+\circ(1))v_k(y),
\quad |y|\le \frac 12M_k^{\frac{2}{n-2}},\quad
\nabla v_k(0)=0.
\end{array}
\right.
\label{ab1}
\end{equation}
Note that (\ref{5-0}) is the same as
\begin{equation}
\max_{ |y|\le \frac 1k M_k^{\frac{2}{n-2}}  }
(v_k(y)|y|^{n-2}) \ge k,
\label{200}
\end{equation}
where $|y|:=\sqrt{ (y^1)^2+\cdots+(y^n)^2}$.
Since we eventually draw contradiction for
large $k$, so throughout the paper $k$ is
large unless otherwise stated.

By standard elliptic estimates,
solutions $v_k$ of
(\ref{ab1}),
 after  passing to a subsequence (still denoted as $v_k$, etc.),
  converge in $C^2_{loc}(\Bbb R^n)$ to
some positive function $U$ satisfying
$U(0)=1, \nabla U(0)=0$ and
$$
-\Delta U=n(n-2) U^{ \frac {n+2}{n-2} }\qquad\mbox{in}\ \Bbb R^n.
$$
By the Liouville type theorem 
 in \cite{CGS},
$$U(y)=(1+|y|^2)^{-\frac {n-2}2}\qquad \mbox{in}\ \Bbb R^n.
$$

For some universal constant $\delta_1>0$, the Green$'$s function
$G(0,x)$ of $-L_g$ on $B(0, 3\delta_1)$,
with respect to zero Dirichlet boundary condition,
is positive and satisfies
$$
\frac 1C dist_g(0,x)^{2-n}\le G(0,x)\le C
dist_g(0,x)^{2-n}, \qquad
x\in B(0, 2\delta_1)\setminus\{0\},
$$
\begin{equation}
\lim_{x\to 0} G(0,x)dist_g(0,x)^{n-2}=\frac 1{  (n-2) \omega_n},
\label{6new}
\end{equation}
where $\omega_n$ denotes the volume of the
standard $(n-1)-$sphere and
  $C>0$ is universal.  By the convergence of $v_k$ to $U$
and the above behavior of the Green$'$s function,
$$
u_k(x)\ge \frac 1C M_k^{-1} G(0,x),
\qquad \mbox{on}\
\partial \left(
B(0,3\delta_1)\setminus B(0, M_k^{ -\frac 2{n-2} })\right).
$$
Since $u_k$ is a supersolution we obtain, using the maximum principle,
$$
u_k(x)\ge C^{-1}M_k^{-1}G(0,x),
\qquad \mbox{on}\
B(0,\delta_1)\setminus B(0, M_k^{ -\frac 2{n-2} }),
$$
i.e.
\begin{equation}
v_k(y)\ge \frac 1{ C (1+|y|^{n-2}) }, \qquad\forall\
0< |y|\le \delta_1 M_k^{\frac{2}{n-2}}.
\label{59}
\end{equation}

For  $\lambda>0$ and for any function $v$, let
$$v^{\lambda}(y):=(\frac{\lambda}{|y|})^{n-2}
v(y^\lambda),  \qquad y^\lambda:=\frac{\lambda^2y}{|y|^2},$$
denote the Kelvin transformation of $v$, and let
$$
\Sigma_{\lambda}:=B(0, \frac 12 M_k^{\frac{2}{n-2}})\setminus
\overline{B(0, \lambda)}=\{y\ |\ \lambda<|y|<
 \frac 12  M_k^{\frac{2}{n-2}}\},
$$
$$
w_{\lambda}(y):= v_k(y)-v_k^{\lambda}(y),\qquad y\in \Sigma_\lambda.
$$

A calculation yields (see \cite{LZ1})
\begin{equation}
\Delta w_{\lambda}+\bar b_i\partial_i w_{\lambda}
+\bar d_{ij}\partial_{ij}w_{\lambda}
-\bar cw_{\lambda}+n(n+2)\xi^{\frac 4{n-2}}w_{\lambda}=
E_{\lambda}, \qquad\mbox{in}\ \Sigma_\lambda,
\label{3eq5}
\end{equation}
where $\xi>0$ is given by
\begin{equation}
n(n+2)\xi^{ \frac 4{n-2}}=\left\{\begin{array}{ll}
n(n-2) \frac { v_k^{ \frac{n+2}{n-2} }-(v_k^\lambda)^{ \frac{n+2}{n-2} } }
{ v_k-v_k^\lambda },\quad v_k\neq v_k^{\lambda},\\
\\
n(n+2)v_k^{\frac{4}{n-2}},\quad v_k=v_k^{\lambda},
\end{array}
\right.
\label{stays}
\end{equation}
and
\begin{eqnarray}
E_{\lambda}&=& \left(\bar c(y)v_k^{\lambda}(y)- (\frac{\lambda}{|y|})^{n+2}
\bar c(y^{\lambda})v_k(y^{\lambda})\right)
-(\bar b_i\partial_i v_k^{\lambda}+\bar d_{ij}\partial_{ij}v^{\lambda}_k)
\nonumber\\
&& +(\frac{\lambda}{|y|})^{n+2}\left(\bar b_i(y^{\lambda})
\partial_i v_k(y^{\lambda})+
\bar d_{ij}(y^{\lambda})\partial_{ij}v_k(y^{\lambda})\right).
\label{eQ}
\end{eqnarray}

By  the convergence of $v_k$ to $U$,
$$
\sigma_k:=
\|v_k-U\|_{C^2(B_2)}\to 0\quad\mbox{as}\ k\to \infty.
$$

\begin{prop}
 For $n\ge 3$, let $v_k$ satisfy
(\ref{ab1}).  Then for
 $\lambda\in (0,2]$ and $y\in \Sigma_\lambda$,
            \begin{equation}
E_\lambda=
 \bar c(y) U^\lambda(y)-
(\frac \lambda {|y|})^{n+2}\bar c(y^\lambda)U(y^\lambda)
+O(1)\sigma_k 
M_k^{ -\frac 4{n-2} } |y|^{-n}.
\label{check}
\end{equation}
where $|O(1)|\le  C_0$ for some 
positive constant $C_0$ independent of $\lambda$, $y$ and $k$.
\label{prop2new}
\end{prop}

\noindent{\bf Proof.}\
For any radially symmetric function $w(y)$,
we have, in conformal normal coordinates,
\begin{equation}
(\Delta_{g_k}-\Delta)w=
(\bar b_i\partial _i+\bar d_{ij}\partial_{ij})w\equiv 0.
\label{rrr}
\end{equation}
Thus
$$
I:=(\bar b_i\partial_i v_k^{\lambda}+\bar d_{ij}\partial_{ij}v^{\lambda}_k)
= 
(\bar b_i\partial_i+\bar d_{ij}\partial_{ij})
[ (v_k-U)^\lambda].
$$

A calculation yields
$$
\partial_i
\bigg\{  (\frac \lambda{|y|})^{n-2}  (v_k-U)(y^\lambda)\bigg\}
=
\partial_i
\bigg\{  (\frac \lambda{|y|})^{n-2} \bigg\}  (v_k-U)(y^\lambda)
+
  (\frac \lambda{|y|})^{n-2}
\partial_i \bigg\{   (v_k-U)(y^\lambda)
\bigg\},
$$
\begin{eqnarray*}
&&
\partial_{ij}
\bigg\{  (\frac \lambda{|y|})^{n-2}  (v_k-U)(y^\lambda)\bigg\}.
\\
&=&
\partial_{ij}
\bigg\{  (\frac \lambda{|y|})^{n-2} \bigg\}  (v_k-U)(y^\lambda)
+
\partial_i
\bigg\{  (\frac \lambda{|y|})^{n-2} \bigg\}
\partial_j
 \bigg\{    (v_k-U)(y^\lambda)  \bigg\}\\
&&
+ \partial_j
\bigg\{  (\frac \lambda{|y|})^{n-2} \bigg\}
\partial_i
 \bigg\{    (v_k-U)(y^\lambda)  \bigg\}
+
(\frac \lambda{|y|})^{n-2}
\partial_{ij}
 \bigg\{    (v_k-U)(y^\lambda)  \bigg\}.
\end{eqnarray*}
It follows, using (\ref{rrr}) and  $\bar d_{ij}\equiv \bar d_{ji}$,
that
\begin{eqnarray*}
I&=&
(\frac \lambda{|y|})^{n-2}
\bar b_i \partial_i
 \bigg\{    (v_k-U)(y^\lambda)  \bigg\}
+2\bar d_{ij} \partial_i
\bigg\{  (\frac \lambda{|y|})^{n-2} \bigg\}
\partial_j
 \bigg\{    (v_k-U)(y^\lambda)  \bigg\}
\\
&&
+
(\frac \lambda{|y|})^{n-2}
\bar d_{ij} \partial_{ij}
 \bigg\{    (v_k-U)(y^\lambda)  \bigg\}.
\end{eqnarray*}
Since $(v_k-U)(0)=0$ and $\nabla (v_k-U)(0)=0$, we have
\begin{equation}
(v_k-U)(y^\lambda)=O(1)\sigma_k|y^\lambda|^2,\qquad
|\nabla (v_k-U)(y^\lambda)|=O(1)\sigma_k|y^\lambda|.
\label{13-1}
\end{equation}

Using (\ref{13-1}) and (\ref{rough}), we obtain
$$
I= O(1)\sigma_k 
M_k^{ -\frac 4{n-2} } |y|^{-n}.
$$
Similarly,
$$
(\frac{\lambda}{|y|})^{n+2}\left(\bar b_i(y^{\lambda})
\partial_i v_k(y^{\lambda})+
\bar d_{ij}(y^{\lambda})\partial_{ij}v_k(y^{\lambda})\right)
= O(1)\sigma_k 
M_k^{ -\frac 4{n-2} } |y|^{-n}.
$$
and
$$
 |\bar c(y)||v_k^\lambda(y)-U^\lambda(y)|
+(\frac \lambda{|y|})^{n+2}
|\bar c(y^\lambda)||v_k(y^\lambda)-U(y^\lambda)|
=O(1)\sigma_k 
M_k^{ -\frac 4{n-2} } |y|^{-n}.
$$
Proposition \ref{prop2new} is established.

\vskip 5pt
\hfill $\Box$
\vskip 5pt

For $\bar l\ge 2$, write the Taylor expansion of $R(x)$ at $0$:
\begin{equation}
R(x)=\sum_{l=2}^{\bar l} \sum_{|\alpha|=l} \frac {\partial_\alpha R}{\alpha!}
x^\alpha+O(|x|^{\bar l+1}).
\label{Rr1}
\end{equation}
Thus, with $r=|y|, y=r\theta$,
\begin{eqnarray}
&& \bar c(y) U^\lambda(y)-
(\frac \lambda {|y|})^{n+2}\bar c(y^\lambda)U(y^\lambda)
\nonumber\\
&=& \sum_{l=2}^{\bar l} M_k^{ -\frac{4+2l}{n-2} }
H_{l,\lambda}(r)
\sum_{|\alpha|=l} \frac {\partial_\alpha R}{\alpha!}
\theta^\alpha
+ O(M_k^{ -\frac 4{n-2} })
|M_k^{ -\frac 2{n-2} }y|^{\bar l+1}
\cdot  (\frac \lambda {|y|})^{n-2},
\label{RR2}
\end{eqnarray}
where
\begin{equation}
H_{l,\lambda}(r)=c(n)\lambda^{n-2} r^{2+l-n}
[1- (\frac \lambda r)^{4+2l}] U(\frac {\lambda^2}r).
\label{RR3}
\end{equation}

Let
\begin{equation}
\bar R^{(l)}:=\frac 1{ |\Bbb S^{n-1}| }
\int_{\theta\in \Bbb S^{n-1} }
\sum_{|\alpha|=l} \frac {\partial_\alpha R}{\alpha!}
\theta^\alpha,
\label{RR6}
\end{equation}
and
\begin{equation}
\tilde  R^{(l)}(\theta):= -\bar R^{(l)}+
\sum_{|\alpha|=l} \frac {\partial_\alpha R}{\alpha!}
\theta^\alpha, \qquad\theta\in \Bbb S^{n-1}.
\label{RR7}
\end{equation}
By (\ref{conformal1}),
\begin{equation}
\bar R^{(2)}=\frac 1{2n}\Delta R=
-\frac 1{12n}|W|^2,
\qquad \mbox{and}\ 
\bar R^{(3)}=0.
\label{SS5}
\end{equation}
We assume that for some constants $\gamma
\ge  0$ and $C\ge 0$,
\begin{equation}
M_k^{ \frac{\gamma}
{(n-2)^2}  } =\circ(M_k^{ \frac 2{n-2}}),\qquad
\|v_k-U\|_{C^2(B(0, 2M_k^{ \frac{\gamma}
{(n-2)^2}  }  ))}\le C M_k^{ -\frac {\gamma}{n-2}}.
\label{RR24}
\end{equation}

We deduce from Proposition \ref{prop2new}, using
(\ref{RR2}), the following
\begin{cor}
For $n\ge 3$,
  let $v_k$ satisfy (\ref{ab1}).
We assume   (\ref{RR24}).
Then, for  $\lambda\in (0, 2]$ and for $y\in \Sigma_\lambda$, we have, for some positive
constant $C_0$ independent of $k$, $\lambda$ and $y$,
\begin{equation}
E_\lambda(y)\le 
 \sum_{l=2}^3 M_k^{ -\frac{4+2l}{n-2} }
H_{l,\lambda}(r)
\tilde R^{(l)}(\theta)
+ C_0 
M_k^{ -\frac {4+\gamma}{n-2} } r^{-n}
+C_0  M_k^{ -\frac{12}{n-2} }
r^{6-n},
\label{RR26}
\end{equation}
\label{cor2-new2}
\end{cor}

In the proof of the following 
 iterated estimates on the rate of 
convergence 
 of $v_k-U$, we make use of some ideas in
 Chen and Lin \cite{ChenLin2} and
 Caffarelli, Hardt and Simon \cite{CHS}, in addition to the
 way of Schoen  \cite{S1} in using the method of moving planes to prove the
 Harnack type inequality (\ref{har}).
\begin{prop} For $n\ge 3$,
 we assume that 
 (\ref{RR24})
 holds  for some
 constants
$0\le\gamma<2(n-2)$,  and $C\ge 0$.   Let $v_k$ 
satisfy
(\ref{ab1}).
Then
  there exist some
positive constants $\delta'>0$ and
  $C_2>0,$
 independent of $k$, such that 
$$
\|v_k-U\|_{C^2(B(0,\delta' R_k)) }
\le C_2 (R_k)^{2-n},
$$
 for any $\{R_k\}$ 
  satisfying, for some $\bar \epsilon\in (0,1)$  
 independent of $k$,
\begin{equation}
2\le R_k =\circ(M_k^{ \frac 2{n-2} }),
\label{SS40}
\end{equation}
\begin{equation}
R_k=\circ(1)
 M_k^{ \frac{4+\gamma}
{  (n-2) (n-2+\bar \epsilon) }
},
\label{SS11}
\end{equation}
\begin{equation}
R_k=
\circ(1)
 M_k^{ \frac {12}{ (n-2)(\max\{6, n-2+\bar\epsilon\})}},
\label{5-1new}
\end{equation}

\begin{equation}
R_k=O(1)
M_k^{ \frac {8}{(n-2)\max\{4+\bar \epsilon, n-2\}} }.
\label{GG8}
\end{equation}
\label{prop5}
\end{prop}
A consequence of Proposition \ref{prop5} is
\begin{cor}
Let $v_k$ satisfy (\ref{ab1}).  For  any $\epsilon>0$,
let
$$
R_k=
\left\{
\begin{array}{rl}
M_k^{ \frac {2-\epsilon}{n-2} },
& \mbox{if}\  3\le n\le 6,\\
M_k^{ \frac 8{  (n-2)^2 }  },
&\mbox{if}\ n\ge 7.
\end{array}
\right.
$$
Then
$$
\limsup_{k\to\infty} (R_k)^{n-2}
\cdot\|v_k-U\|_{ C^2(  B(0, R_k) ) }
<\infty.
$$
\label{cor2-2new}
\end{cor}

\noindent{ \bf Proof of Corollary \ref{cor2-2new}
by using  Proposition \ref{prop5}.}\
Taking first $\gamma = 0$ in (\ref{RR24}) and
 applying Proposition \ref{prop5} with
$R_k=M_k^{ \frac a{  (n-2)^2 }  }$
for $0<a< \min\{4, 2(n-2) \}$,  we see that
 (\ref{RR24}) holds now for any  $0<\gamma<\min\{4, 2(n-2)\}$.
Corollary \ref{cor2-2new} for $n=3,4$ is established. 
For $n\ge 5$, since we now deduce that 
 (\ref{RR24}) holds for any  $0<\gamma<\min\{4, 2(n-2)\}$,
we can apply Proposition \ref{prop5} with
$R_k= M_k^{  \frac a{  (n-2)^2}  }$,
$0<a < \min\{8, 2(n-2)\}$, and $\gamma$ very close to 
$\min\{4, 2(n-2)\}$, and know that
(\ref{RR24}) holds for any $0<\gamma< \min\{8, 2(n-2)\}$.
Corollary \ref{cor2-2new} for $n=5,6$ is established.
For $n\ge 7$, we already know that (\ref{RR24}) holds for any $0<\gamma<8$.
Take $\gamma$ very close to $8$, we can apply
Proposition \ref{prop5} with
$R_k=  M_k^{  \frac 8{  (n-2)^2 }  }$ to 
conclude the proof of  
Corollary \ref{cor2-2new} for $n\ge 7$.

\vskip 5pt
\hfill $\Box$
\vskip 5pt

We prove Proposition \ref{prop5} by the method
of moving spheres, and we need to construct appropriate 
auxiliary
functions to handle the error term $E_\lambda$.
 For $n\ge 3$ and $\alpha<2$,  let
$$
f_{n,\alpha}(r)=
-\frac{1}{(n-\alpha)(2-\alpha)}[r^{2-\alpha}-1]
-\frac 1{(n-\alpha)(n-2)}[r^{2-n}-1],
\quad r\ge 1.
$$
Clearly,
\begin{equation}
f_{n,\alpha}(1)=f_{n,\alpha}'(1)
=0, \qquad 0\le -f_{n,\alpha}(r)
\le C(n,\alpha) r^{2-\alpha},\qquad r\ge 1.
\label{SS1}
\end{equation}
Thinking of $f_{n,\alpha}(r)$ as a radially symmetric function
in $\Bbb R^n$, and let $\Delta$ denote the Laplacian in $\Bbb R^n$, we have
\begin{equation}
\Delta f_{n,\alpha}(r)
=f_{n,\alpha}''(r)+\frac{n-1}{r}f_{n,\alpha}'(r)=-r^{-\alpha},
\quad 
r\ge 1,
\label{21new}
\end{equation}
\begin{equation}
|\frac {d^i}{dr^i} f_{n,\alpha}(r)|
\le C(n,\alpha) |\Delta f_{n,\alpha}(r)| r^{2-i},\qquad
r\ge 1,\  i=0,1,2.
\label{TT3}
\end{equation}

Since
$$
\int_{\theta\in \Bbb S^{n-1} }
\tilde  R^{(l)}(\theta)=0,
$$
we can write
\begin{equation}
\tilde  R^{(l)}(\theta)=\sum_{j=1}^l \sum_{i=1}^{ I_j} a_{ji}^{l}
 Y_j^{(i)}(\theta),
\label{RR9}
\end{equation}
where $Y_j^{(i)}(\theta)$ are spherical harmonics
of degree $j$ satisfying, for some
$\mu_j\ge n-1$,
\begin{equation}
-\Delta_{ \Bbb S^{n-1} } Y_j^{(i)}(\theta)
=\mu_j Y_j^{(i)}(\theta).
\label{RR10}
\end{equation}
Consider, for $\frac 12< \lambda< 2$,
\begin{equation}
\left\{\begin{array}{ll}
\Delta h_{l, j, \lambda}(r)+
(V_{\lambda}(r)-\frac{\mu_j}{r^2})h_{l,j,\lambda}(r)
=-H_{l,\lambda}(r),\ \ \ \lambda<r<
 M_k^{\frac{2}{n-2}},\\
h_{l,j,\lambda}(r)\ge 0,\qquad\qquad \qquad\qquad  \qquad
 \qquad\qquad \ \ \ \lambda<r<
M_k^{\frac{2}{n-2}},\\
h_{l,j,\lambda}(\lambda)=0,
\quad h_{l,j,\lambda}(M_k^{\frac{2}{n-2}})=0,
\end{array}
\right.
\label{eq1aaa}
\end{equation}
where
$$V_{\lambda}(r):=
\left\{
\begin{array}{rll}
n(n-2)
\frac{U(r)^{\frac{n+2}{n-2}}-U^{\lambda}(r)^{\frac{n+2}{n-2}}}
{U(r)-U^{\lambda}(r)}, & & \lambda\ne 1,\\
n(n+2) U(r)^{ \frac 4{n-2} },  & & \lambda=1.
\end{array}
\right.
$$

By Proposition \ref{prop8} in Appendix A,
there  exists some small  $\epsilon_4=
\epsilon_4(n)\in (0,\frac 12)$ such that for $\lambda\in [1-\epsilon_4,
1+\epsilon_4]$,
equation (\ref{eq1aaa}) has a unique classical solution
satisfying
\begin{equation}
\sum_{i=0}^2 |\frac {d^i}{dr^i}
h_{l,j,\lambda}(r)| r^{n-l-4+i}\le C,
 \ \
\lambda<r< M_k^{\frac{2}{n-2}},
\label{eq2aaa}
\end{equation}
where  $C>0$  depends  only  on $n$ and $l$.
From now on we only consider $\lambda$ in this range.
 
Let, with our notation $r=|y|$ and $y=r\theta$,
$$
\tilde h_{l,j,\lambda}^{(i)}(y):=
 M_k^{ -\frac{4+2l}{n-2} }
 h_{l,j,\lambda}(r) a_{ji}^l Y_j^{ (i)}(\theta),
\qquad \lambda\le r\le 4\delta_1 M_k^{ \frac 2{n-2} }.
$$
Then, by (\ref{eq1aaa}),  we have, for
$y\in \Sigma_\lambda$,
\begin{equation}
(\Delta+V_\lambda)\left( \sum_{l=2}^{3}\sum_{j=1}^l
\sum_{i=1}^{I_j}  \tilde h_{l,j,\lambda}^{(i)}(y) \right)
=- \sum_{l=2}^{3}  M_k^{ -\frac{4+2l}{n-2} }
H_{l,\lambda}(r) \tilde R^{(l)}(\theta).
\label{RR22}
\end{equation}

Now we construct the  auxiliary
functions which will be used in the proof of Proposition \ref{prop5}.  
Let, for $y\in \Sigma_\lambda$,
\begin{equation}
\tilde h_{1,\lambda}(y)=
 \sum_{l=2}^{3}\sum_{j=1}^l
\sum_{i=1}^{I_j}  \tilde h_{l,j,\lambda}^{(i)}(y),
\label{RR27}
\end{equation}
\begin{equation}
\tilde h_{3,\lambda}(y)
=Q  M_k^{ -\frac{4+\gamma}{n-2}  }
 f_{n,  2-\hat \epsilon }(\frac r\lambda),   
\label{RR29}
\end{equation}
\begin{equation}
\tilde h_{4,\lambda}(y)
=Q  M_k^{ -\frac{12}{n-2} }
f_{n, \min\{n-6, 2-\hat \epsilon\}}(\frac r\lambda),
\label{RR30}
\end{equation}
and
\begin{equation}
h_\lambda(y)= \tilde h_{1,\lambda}(y)
+\tilde h_{3,\lambda}(y)
+\tilde h_{4,\lambda}(y),
\label{RR31}
\end{equation}
where $\hat \epsilon=  \bar \epsilon/ 9$ and
 $Q>C_0$, independent of $k$, is some large constant
to be fixed later.

Since $M_k^{ -\frac 2{n-2} }|y|=O(1)$
for $y\in \Sigma_\lambda$,
we have, by (\ref{eq2aaa}),
\begin{equation}
|\tilde h_{1,\lambda}(y)|
\le  C\sum_{l=2}^{3} M_k^{ -\frac {4+2l}{n-2} }|h_{l,j,\lambda}(y)|
\le CM_k^{ -\frac 8{n-2} }|y|^{6-n},\quad
y\in \Sigma_\lambda.
\label{SS2}
\end{equation}
Similarly, 
\begin{equation}
|\nabla \tilde h_{1,\lambda}(y)|\le
C M_k^{ -\frac 8{n-2} }|y|^{5-n},
\ \ \ 
|\nabla ^2 \tilde h_{1,\lambda}(y)|\le
C M_k^{ -\frac 8{n-2} }|y|^{4-n}, \qquad y\in \Sigma_\lambda.
\label{SS60}
\end{equation}
By  (\ref{SS1}) ,
we have, for some  $C>0$ independent of $k$ and $Q$,
\begin{equation}
|\tilde h_{3,\lambda}(y)|
\le
CQ M_k^{ -\frac {4+\gamma}{n-2} }
r^{   \hat \epsilon  },
\qquad y\in \Sigma_\lambda,
\label{B2-1}
\end{equation}
\begin{equation}
|\tilde h_{4,\lambda}(y)|
\le  C  Q M_k^{ -\frac{12}{n-2} }
r^{ \max\{8-n, \hat\epsilon\} },
\qquad \forall\ \lambda<|y|<\frac 12 M_k^{ \frac 2{n-2} }.
\label{SS4}
\end{equation}
We also know from  (\ref{21new}) and
 (\ref{TT3}) that, for all $y\in \Sigma_\lambda$,
\begin{equation}
\Delta \tilde h_{3,\lambda}(y)
=- Q
 \lambda^{   -\hat \epsilon  }
M_k^{  -\frac {4+\gamma}{n-2}  }
\cdot r^{ \hat \epsilon-2},
\label{SS21}
\end{equation}
\begin{equation}
\quad \Delta \tilde h_{4,\lambda}(y)
=- Q\lambda^{ \min\{n-8, \hat\epsilon\}  }
M_k^{ -\frac {12}{ n-2} }
r^{ -\min\{n-6, 2-\hat\epsilon\} },
\label{SS22}
\end{equation}
\begin{equation}
|\nabla ^i \tilde h_{m,\lambda}(y)|
\le
 C |y|^{ 2-i}|\Delta \tilde h_{m,\lambda}(y)|,
\qquad   i=0,1,2, \ \ m=3,4.
\label{TT5}
\end{equation}
\begin{lem}  For $n\ge 3$,
 we assume that 
 (\ref{RR24})
 holds for some
 constants 
$0\le\gamma<2(n-2)$  and $C\ge 0$.  
 Let $v_k$
satisfy
(\ref{ab1}), and let  $\{R_k\}$ satisfy (\ref{SS40}),
 (\ref{SS11}) and  (\ref{5-1new}).
  Then
 for  any $\epsilon>0$, there exists $k_0>1$
(can depend on $\epsilon$ and $Q$),
such that for all
$k\ge k_0$,
\begin{equation}
\min_{|y|=r}v_k(y)\le (1+\epsilon)U(r),\qquad\forall\  0<r\le
 R_k.
\label{eq35}
\end{equation}
\label{prop3}
\end{lem}

\noindent{\bf Proof of Lemma \ref{prop3}:}
We prove it by a contradiction argument.
Suppose (\ref{eq35}) is not true, then there exists
 $\epsilon_0>0$, such that
$$\min_{|y|=r_k}v_k(y)> (1+\epsilon_0)U(r_k),$$
for a sequence of
 $r_k\in (0,    R_k]$.
We have written  $U(r)$ for
$U(y)$, $|y|=r$.

By the convergence
of $v_k$ to $U$, we know  $r_k\to \infty$.
Thus
\begin{equation}
\min_{|y|=r_k}v_k(y)\ge (1+\epsilon_0/2)r_k^{2-n}.
\label{B1}
\end{equation}

Fixing 
a small $\epsilon_4'\in (0, \epsilon_4(n))$,
 independent of $k$, such that
$$
U^{\lambda}(y)\le  (1+\frac{\epsilon_0}8)|y|^{2-n}
\qquad\forall\ 0<\lambda\le 1+\epsilon_4',  \
|y|=r_k.
$$
It follows that
\begin{equation}
v_k^{\lambda}(y)\le (1+\frac{\epsilon_0}4)|y|^{2-n},
\qquad \forall\ 0<\lambda\le 1+\epsilon_4', \
 |y|=r_k.
\label{con}
\end{equation}

We will
 derive a contradiction by
applying the method of moving spheres
to $w_\lambda+h_{\lambda}$ with $1-\epsilon_4'\le \lambda\le
1+\epsilon_4'$.

Let
$$
\hat \Sigma_\lambda:=\{ y\ ; |\ \lambda<|y|<r_k\}.
$$

We know that  $h_{\lambda}=0$ on
$\partial B_{\lambda}$ and,  in view of (\ref{SS40}), 
 (\ref{SS11}),  (\ref{5-1new}),
 (\ref{SS4}) and (\ref{SS2}), 
\begin{equation}
h_{\lambda}(y)=\circ(1)|y|^{2-n}, \qquad y\in \hat \Sigma_\lambda,
\label{Z1}
\end{equation}
where $\circ(1)$ denotes some quantity going to zero
as $k\to\infty$, uniform in $y$.

\noindent{\bf Step 1.}\ For $\lambda_0=1-\epsilon_4'$,
\begin{equation}
w_{\lambda_0}(y)+h_{\lambda_0}(y)\ge 0,\qquad
\forall\ \ y\in \hat \Sigma_{\lambda_0}.
\label{67}
\end{equation}

Since $\lambda_0<1$, there exist
some small positive constant
 $\epsilon_5\le \bar \epsilon_0/10$ and some large 
constant
$R_1>10$ such that
\begin{equation}
U(y)-U^{ \lambda_0}(y)
\ge \epsilon_5(|y|-\lambda_0)|y|^{1-n}, \qquad
|y|>\lambda_0,
\label{C2-2}
\end{equation}
\begin{equation}
U(y)>(1-\frac{\epsilon_5}2)|y|^{2-n}, \qquad |y|=R_1,
\label{C2-3}
\end{equation}
\begin{equation}
U^{\lambda_0}(y)<(1-4\epsilon_5)|y|^{2-n}, \qquad
|y|\ge R_1.
\label{C2-4}
\end{equation}

Since $v_k$ converges  in $C^1$
to $U$  in the
region $\lambda_0\le |y|\le R_1$, $h_{\lambda_0}(y)=0$ for
$|y|=\lambda_0$,  and since
$|h_{\lambda_0}(y)|+
|\nabla h_{\lambda_0}(y)|=\circ(1)$
in the same region and uniform in $y$, we deduce from
(\ref{C2-2}), (\ref{C2-3}) and (\ref{C2-4}) that,
for large $k$ as always,
\begin{equation}
w_{\lambda_0}(y)+h_{\lambda_0}(y)>0,\qquad
\lambda_0<|y|\le R_1,
\label{C3-1}
\end{equation}
\begin{equation}
v_k(y)>(1-\epsilon_5)|y|^{2-n}, \qquad
|y|=R_1,
\label{C3-2}
\end{equation}
\begin{equation}
v_k^{ \lambda_0}(y)\le (1-3\epsilon_5)|y|^{2-n}, \qquad
|y|\ge R_1.
\label{C3-3}
\end{equation}

Let $G(0,x)$ be the Greens function of $-L_g$ on $B(0, 3\delta_1)$ as at the
beginning of this section.
Using the maximum principle, we compare
$u_k$ and $(1-\epsilon_5)(n-2)\sigma_n
M_k^{-1} G(0,x)$ as at the
beginning of this section to obtain, by
(\ref{C3-2})
and (\ref{6new}), for some $\delta_2>0$ independent of $k$,
\begin{equation}
v_k(y)\ge (1-2\epsilon_5)|y|^{2-n},
\qquad R_1\le |y|\le \delta_2 
M_k^{ \frac 2{n-2} }.
\label{C4-1}
\end{equation}

Using (\ref{C4-1}), (\ref{C3-3}) and (\ref{Z1}), we obtain
$$
w_{\lambda_0}(y)+h_{\lambda_0}(y)>0,\qquad
R_1\le |y|\le  r_k.
$$
Step 1 follows from this and (\ref{C3-1}).

\medskip

For $\lambda_1=1+\epsilon_4'$, let
$$
\bar\lambda_k=\sup\{\lambda_0\le \lambda\le \lambda_1\
|\ w_\mu+h_\mu\ge 0\ \mbox{in}\
\hat \Sigma_\mu\
\mbox{for all}\ \lambda_0\le \mu\le \lambda\}.
$$

\noindent{\bf Step 2.}\ $\bar \lambda^k=\lambda_1$.

\medskip

Let
$$
\hat O_{\lambda}:=\{y\in \hat \Sigma_{\lambda}\
|\  v_k(y)<2v_k^{\lambda}(y) \}.$$

It follows from (\ref{Z1})  that, for large
$k$ (the largeness of $k$ may depend on $Q$),
\begin{equation}
w_\lambda+h_\lambda(y)=
v_k(y)-v_k^{\lambda}(y)+h_{\lambda}(y)
\ge v_k^{\lambda}(y)+h_{\lambda}(y)
>0
\ \ \mbox{in}\ \ \hat \Sigma_{\lambda}\setminus \hat O_{\lambda},
\label{eq18}
\end{equation}

We also know from (\ref{B1}), (\ref{con}) and (\ref{Z1})
that
\begin{equation}
w_{\bar \lambda^k}(y)
+h_{\bar \lambda^k}(y)>0,\qquad |y|=r_k.
\label{TT50}
\end{equation}

Recall that $v_k$ satisfies (\ref{ab1}),
and $w_\lambda:=v_k-v_k^\lambda$ satisfies
(\ref{3eq5}), with $\xi$ given by (\ref{stays}) and
$E_\lambda$, defined by (\ref{eQ}), satisfying
(\ref{RR26}).
To complete the proof of Lemma \ref{prop3} we need the following

\begin{lem}  For $n\ge 3$,
 we assume that  
 (\ref{RR24})
  holds  for
  some constants 
$0\le\gamma<2(n-2)$  and $C\ge 0$. 
 Let $v_k$
satisfy
(\ref{ab1}) and  let $\{r_k\}=\circ(M_k^{ \frac 2{n-2} })$.
  Then
 we
have,
 for $1-\epsilon_4'\le \lambda\le 1+\epsilon_4'$, and for a large constant $Q$,
\begin{equation}
(\Delta_{g_k}-\bar c+n(n+2)\xi^{\frac{4}{n-2}})h_{\lambda}+E_{\lambda}\le 0
\quad \mbox{in}\quad \hat O_{\lambda}.
\label{eq19}
\end{equation}
\label{prop8-1}
\end{lem}

By (\ref{SS1}) and (\ref{rrr}),
$$
\tilde h_{m,\lambda}\le 0,
\quad \Delta_{g_k} \tilde h_{m,\lambda}\equiv \Delta 
\tilde h_{m,\lambda},\qquad m=3,4.
$$
It follows, using 
(\ref{SS21}),
(\ref{SS22}) and the smallness of
 $|\lambda-1|\le \epsilon_4'$,
 that, for $y\in \Sigma_\lambda$,
\begin{equation}
 (\Delta_{g_k}+n(n+2)\xi^{ \frac 4{n-2}})
\tilde h_{3, \lambda}
\le -\frac Q2
 M_k^{   -\frac{4+\gamma}{n-2}  }
\cdot r^{ \hat \epsilon-2},
\label{SS27}
\end{equation}
\begin{equation}
 (\Delta_{g_k}+n(n+2)\xi^{ \frac 4{n-2} })
\tilde h_{4, \lambda}
\le  -\frac Q2
 M_k^{ -\frac {12}{ n-2} }
r^{ -\min\{n-6, 2-\hat\epsilon\}}.
\label{SS28}
\end{equation}

Using (\ref{rough}), (\ref{SS40}) and (\ref{TT5}),
we have, for $y\in \hat \Sigma_\lambda$,
\begin{equation}
|\bar c\tilde h_{m,\lambda}(y)|
\le CM_k^{ -\frac 8{n-2} }|y|^2
\cdot |y|^{2}|\Delta \tilde  h_{m,\lambda}(y)|
=\circ(1) |\Delta \tilde  h_{m,\lambda}(y)|,
\ \ m=3,4.
\label{TT12}
\end{equation}
Putting together the above four estimates,
we have
\begin{equation}
\sum_{m=3}^4
(\Delta_{g_k}-\bar c+n(n+2)\xi^{ \frac 4{n-2}})
\tilde h_{m, \lambda}
\le 
 -\frac Q4 D_{k, \hat \epsilon}(r),
\qquad y\in \hat \Sigma_\lambda.
\label{SS30}
\end{equation}
where 
$$
D_{k, \hat \epsilon}(r):= M_k^{   -\frac{4+\gamma}{n-2}  }
\cdot r^{ \hat \epsilon-2}
+ M_k^{ -\frac {12}{ n-2} }
r^{ -\min\{n-6, 2-\hat\epsilon\}}.
$$
The right hand side of (\ref{SS30}) is
of good sign and  will be used to absorb
other terms.

Using (\ref{rough}), (\ref{SS2})
and (\ref{SS60}),
we have,
for $y\in \hat \Sigma_\lambda$,
\begin{equation}
|\bar c||\tilde h_{1,\lambda}(y)|
+|\bar b_i \partial_i \tilde h_{1,\lambda}(y)|
+|\bar d_{ij}\partial_{ij}  \tilde h_{1,\lambda}(y)|
\le C  M_k^{-\frac {12}{n-2} }
|y|^{6-n}\le C D_{k, \hat \epsilon}(r).
\label{TT33}
\end{equation}

We  give an estimate of $\xi$, given by (\ref{stays}), in the following
\begin{lem} 
For  $n\ge 3$, we assume (\ref{RR24}) for some $0\le \gamma
\le 2(n-2)$.   Let $v_k$
satisfy
(\ref{ab1}).  Then, 
there exists $C$, independent of $k$, such that
\begin{equation}
|n(n+2)\xi^{\frac{4}{n-2}}(y)-V_{\lambda}(|y|)|\le
CM_k^{-\frac{\gamma}{n-2}}
|y|^{n-6},\qquad
\lambda\le |y|\le 2M_k^{  \frac{\gamma}{(n-2)^2}  },
\label{est1}
\end{equation}
and
\begin{equation}
|n(n+2)\xi^{\frac{4}{n-2}}-V_{\lambda}|\le C|y|^{-4},\quad
y\in \hat O_{\lambda}.
\label{est2}
\end{equation}
\label{lem8-1}
\end{lem}

\noindent{\bf Proof of Lemma \ref{lem8-1}.}\
By (\ref{RR24}),
$$v_k(y)=U(y)+a(y),\qquad v_k^{\lambda}(y)=U^{\lambda}(y)+b(y),
\qquad \lambda\le |y|\le 2M_k^{  \frac{\gamma}{(n-2)^2}  },
$$
where $a(y)$ and $b(y)$ satisfy
$$
|a(y)|+|b(y)|\le CM_k^{-\frac{\gamma}{n-2}},
\qquad  \lambda\le |y|\le 2M_k^{  \frac{\gamma}{(n-2)^2}  }.
$$
Then, with $n^*=\frac{n+2}{n-2}$ and for
$\lambda\le |y|\le 2M_k^{ \frac \gamma{ (n-2)^2} }$,
\begin{eqnarray*}
\frac{v_k^{n^*}-v_k^{\lambda}(y)^{n^*}}{v_k-v_k^{\lambda}}
&=& \frac{(U+a)^{n^*}-(U^{\lambda}+b)^{n^*}}{(U+a)-(U^{\lambda}+b)}\\
&=&\frac{\int_0^1\frac{d}{dt}\{(t(U+a)+(1-t)(U^{\lambda}+b))^{n^*}\}dt}
{(U+a)-(U^{\lambda}+b)}\\
&=&n^*\int_0^1(t(U+a)+(1-t)(U^{\lambda}+b))^{\frac{4}{n-2}}dt \\
&=&n^*\int_0^1(tU+(1-t)U^{\lambda})^{\frac{4}{n-2}}dt
+O(1)(|a(y)|+|b(y)|)|y|^{n-6}\\
&=& \frac 1{ n(n-2)} V_\lambda+O(M_k^{-\frac{\gamma}{n-2}}|y|^{n-6}).
\end{eqnarray*}
Estimate (\ref{est1}) is established.
Estimate (\ref{est2}) is obvious.

\vskip 5pt
\hfill $\Box$
\vskip 5pt

By (\ref{est1}) and (\ref{SS2}), we have,
for $\lambda\le |y|\le 2M_k^{ \frac{\gamma}{ (n-2)^2 }}$,
\begin{equation}
|n(n+2)\xi^{ \frac 4{n-2} }-V_\lambda||\tilde h_{1,\lambda}|
\le C  M_k^{ -\frac {8+\gamma}{n-2} }\le C  M_k^{ -\frac {4+\gamma}{n-2} }
r^{-2}
=O(1) D_{k, \hat \epsilon}(r).
\label{TT16}
\end{equation}
By (\ref{est2}),
(\ref{SS2}) and the fact $\gamma\le 2(n-2)$,
we have, for
$y\in
 \hat O_\lambda$ and $|y|\ge  2M_k^{ \frac{\gamma}{ (n-2)^2 }}$,
\begin{equation}
|n(n+2)\xi^{ \frac 4{n-2} }-V_\lambda||\tilde h_{1,\lambda}|
\le C  M_k^{ -\frac 8 {n-2} } |y|^{2-n}
\le C  M_k^{ -\frac {4+\gamma} {n-2} } |y|^{-2}
=O(1) D_{k, \hat \epsilon}(r).
\label{TT17}
\end{equation}
By  (\ref{RR27}),  (\ref{RR22}),  (\ref{TT16})
and (\ref{TT17}),
\begin{eqnarray*}
&&
(\Delta+n(n+2)\xi^{\frac 4{n-2} })\tilde h_{1,\lambda}
\\
&\le & (\Delta+V_\lambda)\tilde h_{1,\lambda}
+|n(n+2)\xi^{\frac 4{n-2} }-V_\lambda||\tilde h_{1,\lambda}|
\\
&\le &
- \sum_{l=2}^{3}  M_k^{ -\frac{4+2l}{n-2} }
H_{l,\lambda}(r) \tilde R^{(l)}(\theta)
+ O(1) D_{k, \hat \epsilon}(r),
\qquad \mbox{in}\  \hat O_\lambda.
\end{eqnarray*}
Thus, in view of (\ref{TT33}),
\begin{eqnarray}
&&
(\Delta_{g_k}-\bar c+n(n+2)\xi^{\frac 4{n-2} })\tilde h_{1,\lambda}
\nonumber\\
&\le& - \sum_{l=2}^{3}  M_k^{ -\frac{4+2l}{n-2} }
H_{l,\lambda}(r) \tilde R^{(l)}(\theta)
+C D_{k, \hat \epsilon}(r),
\qquad
\mbox{in}\  \hat O_\lambda.
\label{TT35}
\end{eqnarray}

Then we fix  a large constant $Q$, estimate (\ref{eq19}) follows from
(\ref{RR26}), (\ref{TT35}) and (\ref{SS30}).  Lemma \ref{prop8-1}
is established.

\vskip 5pt
\hfill $\Box$
\vskip 5pt

Now we establish Step 2:  We know that
$w_{\bar \lambda^k}+h_{\bar \lambda^k}$ is non-negative
in $ \hat \Sigma_{\bar \lambda^k}$,
and, by (\ref{eq19}) ,  satisfies
\begin{equation}
\big(\Delta_{g_k}-\bar c+n(n+2)\xi^{ \frac 4{n-2} }\big)
(w_{\bar \lambda^k}+h_{\bar \lambda^k})\le 0, \qquad
\mbox{in}\ \hat O_{\bar \lambda^k}.
\label{ehh3}
\end{equation}
Since  $w_{\bar \lambda^k}+h_{\bar \lambda^k}$
satisfies
(\ref{TT50}) and  (\ref{eq18}) with $\lambda=\bar \lambda^k$,
we apply the strong  maximum principle and the Hopf lemma   to obtain
$$
w_{\bar \lambda^k}+h_{\bar \lambda^k}>0\qquad
\mbox{in}\ \hat  \Sigma_{ \bar\lambda^k},
$$
and
$$
\frac {\partial }{ \partial \nu}(w_{\bar \lambda^k}+h_{\bar \lambda^k})
>0\qquad\mbox{on}\ \partial B(0, \bar \lambda^k),
$$
where $\frac {\partial }{ \partial \nu}$ denotes the differentiation
in the outer normal direction.
In view of (\ref{TT50}) and  the above two
 estimates above, we must have
$\bar \lambda^k=\lambda_1$.
Step 2 is established.

\medskip

By Step 2,
$w_{\lambda_1}+h_{\lambda_1}\ge 0$
in $\hat \Sigma_{\lambda_1}$.  Sending
$k$ to infinity, we obtain
\begin{equation}
U^{\lambda_1}(y)\le U(y), \qquad \forall \ |y|\ge \lambda_1.
\label{U11}
\end{equation}
But the above is not satisfied by $U$, a fact easily checked
using $\lambda_1>1$.
This leads to contradiction. 
Lemma \ref{prop3} is established.

\vskip 5pt
\hfill $\Box$
\vskip 5pt

\begin{lem} Under the hypotheses of Lemma \ref{prop3} ,
there exist $\delta\in (0,1)$ and $C>0$, independent of $k$,
 such that for all large $k$,
$$v_k(y)\le CU(y),\qquad |y|\le \delta R_k.
$$
\label{prop4}
\end{lem}

A consequence of Proposition \ref{prop5} and Lemma \ref{prop4} is
\begin{cor} For $n\ge 8$ and for any $\epsilon>0$,
 there exists some  constant
$C>1$, independent of $k$,
such that
$$
v_k(y)\le C U(y),
\qquad \forall \ |y|\le M_k^{ \frac{12-\epsilon}{  (n-2)^2 }  }.
$$
\label{cor2-5}
\end{cor}

\noindent{\bf Proof of Corollary \ref{cor2-5}.}\
By Corollary \ref{cor2-2new}, (\ref{RR24}) is satisfied
with $\gamma=8$.
For any $\epsilon>0$, let $R_k= M_k^{ \frac{12-\epsilon}{  (n-2)^2 }  }$,
and let $\bar \epsilon>0$ be
sufficiently small 
(depending on $\epsilon$).  Using $n\ge 8$, we easily
check that hypotheses of Lemma \ref{prop4} are satisfied,
and
Corollary \ref{cor2-5} follows from the lemma.

\vskip 5pt
\hfill $\Box$
\vskip 5pt

\noindent{\bf Proof of Lemma \ref{prop4}:}  The proof is very similar to
the proof of lemma 3.2 in \cite{ChenLin2}.
Let $G_k$ (will be denoted as $G$)
 be the Green's function
of $-L_{g_k}$ on $B(0, \delta_1 M_k^{ \frac 2{n-2}})$
with respect to zero Dirichlet boundary data, where $\delta_1$
is the constant above (\ref{6new}),
and let $y_1$ be a minimum
point of  $v_k$  on $|y|=R_k$.
For $\epsilon>0$,
there exists some constant $\delta\in (0,1)$,
independent of $k$,
such that, for large $k$, the following estimates
go through:
\begin{eqnarray*}
v_k(y_1)&\ge &\int_{B(0,\delta_1 M_k^{ \frac 2{n-2}})}
G(y_1,\eta)n(n-2)
v_k(\eta)^{\frac{n+2}{n-2}}dV_{g_k}
\\
&\ge& \int_{B(0,\delta R_k)}G(y_1,\eta)n(n-2)
v_k(\eta)^{\frac{n+2}{n-2}}dV_{g_k},
\end{eqnarray*}
 and, using (\ref{SS40}),
$$G(y_1,\eta)\ge \frac{(1-\epsilon/2)}{(n-2)\omega_n}|y_1-\eta |^{2-n}
\ge \frac{(1-3\epsilon/4)}{(n-2)\omega_n}|y_1|^{2-n},
\qquad |\eta|= \delta R_k, $$
where $\omega_n$ denotes
the volume of the standard $(n-1)-$sphere.
Since $dV_{g_k}=(1+\circ(1))d\eta$, we have
$$v_k(y_1)\ge \frac{(1-\epsilon)n}{\omega_n}|y_1|^{2-n}
\int_{B(0,\delta R_k)}v_k^{\frac{n+2}{n-2}}(\eta)d\eta $$
On the other hand, by Lemma \ref{prop3},
$$
v_k(y_1)\le  (1+\epsilon) U(y_1)\le (1+2\epsilon)|y_1|^{2-n}.
$$
So
$$\int_{B(0,\delta R_k)}v_k^{\frac{n+2}{n-2}}(\eta)d\eta\le
(1+4\epsilon)\omega_n/n.$$
A direct computation gives,
$$\int_{\Bbb R^n}U^{\frac{n+2}{n-2}}=\frac{ \omega_n}{ n }.
$$
By the convergence of $v_k$ to $U$,
there exists  some $R_1$, depending only on $n$ and $\epsilon$,
 such that, for large $k$,
$$\int_{R_1\le |\eta|\le \delta R_k}v_k^{\frac{n+2}{n-2}}d\eta \le 5\epsilon.
$$
Using the second line of  (\ref{ab1}),
$$\int_{R_1\le |\eta|\le \delta R_k}v_k^{\frac{2n}{n-2}}d\eta \le
(\bar b+1)\int_{R_1\le |\eta|\le \delta R_k}v_k^{\frac{n+2}{n-2}}d\eta
\le 5(\bar b+1)\epsilon.$$

For each $2R_1<r<\delta R_k/2$, we consider
$\tilde v_k(z)=r^{\frac{n-2}2}v_k(rz)$ for $1/2<|z|<2$.
Then  $\tilde v_k$ satisfies
$$\frac{1}{\sqrt{g(rz)}}\partial_{z_i}(\sqrt{g(rz)}g^{ij}(rz)\partial_{z_j}
\tilde v_k(z))-\bar cr^2\tilde v_k(z)+n(n-2)\tilde v_k(z)^{\frac{n+2}{n-2}}
=0,\quad 1/2<|z|<2.$$
We know that
$\int_{ \frac 12\le |z|\le 2}
\tilde v_k(z)^{ \frac {2n}{n-2} }\le 5(\bar b +1)\epsilon$.
Fix some universally small  $\epsilon>0$, we apply
the  Moser iteration technique to obtain
$\tilde v_k(z)\le C$ for $\frac 34\le |z|\le \frac 43$,
where $C$ is independent of $k$.
With this, we apply the Harnack inequality to obtain
$\max_{ |z|=1} \tilde v_k(z)
\le C\min_{ |z|=1} \tilde v_k(z)$, i.e.,
$\max_{ |y|=r} v_k(y) \le C\min_{ |y|=r} v_k(y)$.
By Lemma \ref{prop3},
$$\min_{ |y|=r} v_k(y) \le (1+\epsilon)U(r).$$
Lemma \ref{prop4} follows from these
together with the convergence of $v_k$ to $U$.

\vskip 5pt
\hfill $\Box$
\vskip 5pt

\noindent{\bf Proof of Proposition \ref{prop5}:}
The argument below is very similar to the proof of lemma 3.3
in \cite{ChenLin2}.
Let $
\displaystyle{
\Lambda_k=\max_{ |y|\le \delta R_k}|(v_k-U)(y)|
}
$ 
  and let
$w_k=\Lambda_k^{-1}(v_k-U)$. We will  show that
\begin{equation}
\Lambda_k\le C_2 R_k^{2-n}
\label{eq31}
\end{equation}
 for some $C_2>0$, independent of $k$.
Suppose this is not true, then, along a subsequence,
\begin{equation}
\Lambda_kR_k^{n-2}\to \infty.
\label{cc1}
\end{equation}
By Lemma \ref{prop4},
\begin{equation}
w_k(y)\le C\Lambda_k^{-1} U(y), \qquad |y|\le \delta R_k.
\label{going0}
\end{equation}
By (\ref{cc1}) and (\ref{going0}),
\begin{equation}
\max_{ \partial B(0, \delta R_k) }|w_k|\to 0.
\label{going}
\end{equation}
Since
\begin{eqnarray*}
\Delta_{g_k}v_k-\bar cv_k+n(n-2)v_k^{\frac{n+2}{n-2}}=0,\qquad \quad
|y|\le \delta_1 M_k^{\frac{2}{n-2}},\\
\Delta_{g_k}U-\bar cU+n(n-2)U^{\frac{n+2}{n-2}}=
-\bar cU, \quad |y|\le
\delta_1 M_k^{\frac{2}{n-2}},
\end{eqnarray*}
 $w_k$ satisfies
$$\Delta_{g_k}w_k-\bar cw_k=
-n(n+2)\hat \xi^{\frac{4}{n-2}}w_k+
\Lambda_k^{-1}\bar cU,
\quad \mbox{in}\quad |y|\le \delta R_k $$
where $\hat \xi$ is between $v_k$ and $U$.
By Lemma \ref{prop4}, there is $C>0$ such that
$$|\hat \xi(y)|\le C(1+|y|)^{2-n},\qquad |y|\le \delta R_k.$$
By (\ref{cc1}),
$\Lambda_k^{-1}=\circ(1)R_k^{n-2}$.

Fixing $\epsilon>0$ sufficiently small, 
we
 know from (\ref{rough}) and
(\ref{GG8})
 that, for $|y|\le \delta R_k$,
$$
|\bar c(y)| R_k^{ 2+\epsilon}+|\bar b_i(y)|  R_k^{ 1+\epsilon}+
|\bar d_{ij}(y)| R_k^{\epsilon}
\le CM_k^{ -\frac 4{n-2} }  R_k^{ 2+\epsilon}
=\circ(1).
$$

By (\ref{cc1}), $\Lambda_k^{-1}=\circ(1)R_k^{n-2}$.
So we have, using (\ref{rough}) and (\ref{GG8}),
\begin{equation}
\Lambda_k^{-1} |\bar cU|
\le \circ(1) 
M_k^{ -\frac {8}{n-2} }
R_k^{ \max\{4+\epsilon, n-2\} }  (1+|y|)^{-2-\epsilon}
=\circ(1) (1+|y|)^{-2-\epsilon}.
\label{GG9} 
\end{equation}

It follows that $w_k$ satisfies, for $|y|\le \delta R_k$, that 
$$
\left( \Delta+
\frac{\circ(1) \partial_{ij} }{(1+|y|)^\epsilon}
+\frac { \circ(1) \partial_i  }{ (1+|y|)^{1+\epsilon} }
+\frac  { \circ(1) }{  (1+|y|)^{2+\epsilon} }\right) w_k(y)
=O(1)(1+|y|)^{ -2-\epsilon}.
$$
Let $\eta(r)=(1+r^2)^{ -\frac\epsilon 2}$, we have,
for some $C>0$ depending only on $\epsilon$ and $n$,
$$
\left( \Delta+
\frac{\circ(1) \partial_{ij} }{(1+|y|)^\epsilon}
+\frac { \circ(1)  \partial_i }{ (1+|y|)^{1+\epsilon} }
+\frac  { \circ(1)   }{  (1+|y|)^{2+\epsilon} }\right) 
\eta\le -C^{-1}(1+r)^{ -\epsilon-2}.
$$

Taking a large constant positive $Q$, independent of $k$,
we have, for $|y|\le \delta R_k$, 
$$
\left( \Delta+
\frac{\circ(1) \partial_{ij}}{(1+|y|)^\epsilon}
+\frac { \circ(1) \partial_i }{ (1+|y|)^{1+\epsilon} }
+\frac  { \circ(1) }{  (1+|y|)^{2+\epsilon} }\right) 
(\pm w_k- \max_{ |z|=\delta R_k }|w_k(z)|-Q\eta)\ge 0.
$$
By  the maximum principle,
$$
|w_k(y)|\le \eta (|y|)+ \max_{ |z|=\delta R_k }|w_k(z)|,
\qquad |y|\le \delta R_k.$$

Next, we deduce, using standard elliptic estimates,
from  (\ref{GG9}) and the equation of $w_k$   
 that $w_k$ converges
in $C^2_{loc}(\Bbb R^n)$ to some $w_0$ satisfying
$$
\left\{
\begin{array}{l}
 \Delta w_0+n(n+2)U^{\frac 4{n-2}}w_0=0,\quad \mbox{in}\ \Bbb R^n,
\\
w_0(0)=0,\ \nabla w_0(0)=0, \  \lim_{|y|\to \infty}w_0(y)=0.
\end{array}
\right.
$$
By lemma 2.4 in \cite{ChenLin2}, 
$w_0\equiv 0$.

Let  $y_k$ be a maximum
point of $w_k(y)$ in $|y|\le \delta R_k$, i.e.,
$w_k(y_k)=1, |y_k|\le \delta R_k$.
By the above estimates,
$$1=w_k(y_k)\le C(1+|y_k|)^{-1}+\circ (1),
$$
so, $\{y_k\}$ must be bounded, and therefore, by the convergence
of $w_k$ to $w_0\equiv 0$, $w_k(y_k)\to 0$.  This contradicts to
$w_k(y_k)=1$.  Thus we have established
(\ref{eq31}).
Using the equation satisfied by $w_k$, we have, by
standard elliptic theories,
$$
\|w_k\|_{ C^2(B(0,\delta R_k-1))}\le C
\|w_k\|_{ L^{\infty}(B(0,\delta R_k))}\le C.
$$
Proposition \ref{prop5} follows from this in view of
(\ref{eq31}).

\vskip 5pt
\hfill $\Box$
\vskip 5pt

Now we give the  

\noindent{\bf Proof of 
 (\ref{eq4})  
 for $3\le n\le 7$.}\  For this, we only need to 
reach a contradiction to (\ref{5-0}) in dimension
$3\le n\le 7$. 
As pointed before,
(\ref{5-0}) is equivalent to  (\ref{200}).
It is easy to see from (\ref{200}) that
 $\frac 1k M_k^{ \frac 2{n-2} }\to\infty$.
We know from Corollary \ref{cor2-5} that 
(\ref{RR24}) holds in dimensions $3\le n\le 6$
for any $0<\gamma <2(n-2)$ while in dimension $n=7$ it holds
for $\gamma=8$.  Let,
in dimensions $3\le n\le 7$, $R_k= k^{-\frac 14}
M_k^{\frac 2{n-2}}$.
Then $\{R_k\}$ satisfy  (\ref{SS40}),
 (\ref{SS11}) and (\ref{5-1new})
 with
the above $\gamma$ and sufficiently small
$\bar \epsilon$.  Thus, by Lemma \ref{prop4},
$$
v_k(y)\le CU(y),\qquad
\forall\ |y|\le
k^{-\frac 12}  M_k^{ \frac 2{n-2} },
$$
where $C>0$ is independent of $k$.  This violates
(\ref{200}).  Thus
estimate (\ref{eq4}) in dimension 
$3\le n\le 7$ is established.

\vskip 5pt
\hfill $\Box$
\vskip 5pt

The following is a  Pohozaev type identity.
\begin{lem}
For $n\ge 3$,
let $u$ be a solution of the  Yamabe equation, then in a neighborhood of any
point $P\in M$, the following identity holds in a normal coordinate of $P$.
\begin{eqnarray*}
&&\int_{|x|\le \sigma}\{(-b_i\partial_iu-d_{ij}\partial_{ij}u)(\nabla u\cdot x
+\frac{n-2}2u)
-\frac{c(n)}2u^2(x\cdot \nabla R(x))
-c(n)R(x)u^2\}\\
&&+\frac{\sigma}2c(n)\int_{|x|=\sigma}R(x)u^2
-\frac{(n-2)^2}2\sigma\int_{|x|=\sigma}u^{\frac{2n}{n-2}}
=B(\sigma,u,\nabla u)
\end{eqnarray*}
where
$$B(\sigma,u,\nabla u)=\int_{|x|=\sigma}
(|\frac{\partial u}{\partial \nu}|^2\sigma-\frac 12|\nabla u|^2\sigma
+\frac{n-2}2u\frac{\partial u}{\partial \nu}) $$
\label{lemcom1}
\end{lem}

\noindent{\bf Proof.}\ This 
identity is established for $n=3$
in \cite{LiZhu}.  A modification of the proof there yields the above lemma.

\vskip 5pt
\hfill $\Box$
\vskip 5pt

Applying Lemma \ref{lemcom1} to $u=u_k$
(see (\ref{vk})) with
$\sigma=M_k^{ -\frac 2{n-2} }R_k'$, we have,
after a change of variables,
\begin{equation}
I_1[v_k]+I_2[v_k]+I_3[v_k]+I_4[v_k]=I_5[v_k],
\label{pohov}
\end{equation}
$$
I_1[v_k]=\int_{|y|\le R_k'}
(-\bar b_i \partial_i v_k-\bar d_{ij}\partial_{ij}v_k)
(\nabla v_k\cdot y+\frac {n-2}2 v_k),
$$
$$
I_2[v_k]=-\frac {c(n)}2
M_k^{ -\frac 4{n-2} }
\int_{|y|\le R_k'}
\bigg\{ (M_k^{ -\frac 2{n-2} }y)\cdot \nabla R(M_k^{ -\frac 2{n-2} }y)
+2R(M_k^{ -\frac 2{n-2} }y)\bigg\}
v_k^2(y),
$$
$$
I_3[v_k]=  \frac {c(n)}2
M_k^{ -\frac 4{n-2} }   R_k' 
 \int_{  |y|=R_k' }
R(M_k^{ -\frac 2{n-2} }y)v_k^2(y),
$$
$$
I_4[v_k]= -\frac{  (n-2)^2}2 R_k'
 \int_{  |y|=R_k' }
v_k(y)^{  \frac {2n}{n-2}},
$$
$$
I_5[v_k]=
\int_{ |y|=R_k' }
\bigg\{  (|\frac{\partial v_k}{ \partial \nu}|^2
-\frac 12|\nabla v_k|^2) R_k'
+\frac{n-2}2
v_k \frac{\partial v_k}{ \partial \nu}\bigg\}
= O(1) (R_k')^{2-n}.
$$

Let $\beta_2, \beta_4, \beta_2'', \beta_2''', \beta_3'''\ge 0 $
satisfy,
  for some constant
$C\ge 0$,
\begin{equation}
\bar R^{(2)}\le
C M_k^{ -\frac {\beta_2}{n-2}},\quad
\bar R^{(4)}\le
C M_k^{ -\frac {\beta_4}{n-2}},
\label{RR23}
\end{equation}
\begin{equation}
|b_i(x)|\le C
M_k^{ -\frac {\beta_2''}{n-2} }|x|^2+C|x|^3,
\quad |d_{ij}(x)|\le C \sum_{ l=2}^3
M_k^{ -\frac {\beta_l'''}{n-2} }|x|^{l}+C|x|^4,
\label{TT1}
\end{equation}
or, equivalently,
\begin{equation}
|\bar b_i(y)|\le
C
M_k^{ -\frac {6+\beta_2''}{n-2} }
|y|^2+C M_k^{ -\frac {8}{n-2} }
|y|^3, \quad
|\bar d_{ij}(y)|\le C
 \sum_{ l=2}^3M_k^{ -\frac {2l+\beta_l'''}{n-2} }|y|^l
+C M_k^{ -\frac {8}{n-2} }
|y|^4.
\label{TT2}
\end{equation}
We will always take $\beta_3'''=\beta_2''$.

\begin{lem} For $n\ge 7$, $\bar l\ge 2$,  let
$v_k$ satisfy  the first line of (\ref{ab1}), and we assume
(\ref{TT1}) holds for some constants $\beta_2'', \beta_2''', 
\beta_3'''\ge 0$.
For $2\le R_k'\le \frac 14 M_k^{ \frac 2{n-2}}$,
we assume, for some constants
$\gamma_1, \gamma_2, C\ge 0$,
\begin{equation}
v_k(y)\le CU(y),\qquad
|y|\le 2R_k',
\label{new99}
\end{equation}
and
\begin{equation}
|\nabla^j (v_k-U)|\le C  M_k^{ -\gamma_1 }
(1+|y|)^{-\gamma_2-j},
\qquad |y|\le R_k', j=0,1,2.
\label{1-e1}
\end{equation}
Then
\begin{eqnarray}
&&
-\frac {  c(n)  }2 |\Bbb S^{n-1}|
\sum_{l=2}^{\bar l}  (l+2)
M_k^{ -\frac{4+2l}{n-2} }
\bar R^{ (l)}\int_0^{ R_k'} r^{l+n-1}U(r)^2dr\nonumber\\
&=& I_5[v_k]
+O(1)  \sum_{l=2}^{  \bar l}
M_k^{ -\gamma_1-\frac{4+2l}{n-2} }\int_{ |y|\le R_k' }
(1+|y|)^{ 2-n-\gamma_2 +l}
\nonumber\\
&&+O(1) M_k^{ -\frac{6+2\bar l}{n-2} }\int_{ |y|\le R_k' }
(1+|y|)^{ 5-2n+\bar l}\nonumber\\
&&
+ O(1)M_k^{ -\gamma_1- \frac {6+\beta_2''}{n-2} }
\int_{|y|\le R_k'} (1+|y|)^{ 3-n-\gamma_2}\nonumber\\
&&
+O(1) M_k^{ -\gamma_1-\frac 8{n-2} }
\int_{ |y|\le R_k'}
(1+|y|)^{ 4-n-\gamma_2}\nonumber\\
&& + O(1)\sum_{l=2}^3M_k^{ -\gamma_1- \frac {2l+\beta_l'''}{n-2} }
\int_{|y|\le R_k'} (1+|y|)^{ -n+l-\gamma_2}\nonumber\\
&& + I_3[v_k]
+ O(1)
(R_k')^{-n},
\label{2e-1}
\end{eqnarray}
and
\begin{equation}
I_5[v_k]= O(1) (R_k')^{2-n},\quad I_3[v_k]=O(1)M_k^{-\frac{8}{n-2}}(R_k')^{6-n}.
\label{3-e0}
\end{equation}
\label{cor2-6}
\end{lem}

\noindent{\bf Proof.}\ 
Applying standard elliptic estimates to the equation of
$v_k$ and using (\ref{new99}), we obtain
\begin{equation}
|\nabla v_k(y)|\le C (1+|y|)^{ 1-n},
\quad |\nabla^2 v_k(y)|\le C (1+|y|)^{-n},
\qquad |y|\le R_k',
\label{88}
\end{equation}
where,
 and throughout
 the proof, $C$ denotes various positive constants independent of
$k$. 

The desired estimates are  deduced from
(\ref{pohov}).  The 
main term there is 
\begin{eqnarray}
I_2[v_k]
&=& -\frac {  c(n)}2 
\sum_{l=2}^{  \bar l}
\sum_{ |\alpha|=l}
\int_{ |y|\le R_k' }
\bigg\{ (\frac {(l+2)}{ \alpha!})M_k^{ -\frac{4+2l}{n-2} }
(\partial _\alpha R) y^\alpha
+O(1) M_k^{ -\frac{ 6+2\bar l}{n-2} }|y|^{ \bar l+1}
\bigg\} v_k^2 \nonumber \\
&=& -\frac {  c(n)}2 \sum_{l=2}^{  \bar l}
\sum_{ |\alpha|=l}
 (\frac {(l+2)}{ \alpha!})M_k^{ -\frac{4+2l}{n-2} }
\int_{ |y|\le R_k' }
(\partial _\alpha R) y^\alpha
U(y)^2 \nonumber \\
&&
+O(1)\sum_{l=2}^{  \bar l}
M_k^{ -\gamma_1-\frac{4+2l}{n-2} }\int_{ |y|\le R_k' }
(1+|y|)^{ 2-n-\gamma_2 +l}
\nonumber \\
&&+O(1) M_k^{ -\frac{6+2\bar l}{n-2} }\int_{ |y|\le R_k' }
(1+|y|)^{ 5-2n+\bar l} \nonumber \\
&=&  -\frac {  c(n)}2 \sum_{l=2}^{  \bar l}
(l+2)M_k^{ -\frac{4+2l}{n-2} }
|\Bbb S^{n-1}| \bar R^{ (l)}
\int_0^{ R_k' } r^{l+n-1}U(r)^2dr \nonumber \\
&&+O(1)  \sum_{l=2}^{  \bar l}
M_k^{ -\gamma_1-\frac{4+2l}{n-2} }\int_{ |y|\le R_k' }
(1+|y|)^{ 2-n-\gamma_2 +l}
\nonumber \\
&&+O(1) M_k^{ -\frac{6+2\bar l}{n-2} }\int_{ |y|\le R_k' }
(1+|y|)^{ 5-2n+\bar l}.
\label{i2crude}
\end{eqnarray}

 Using (\ref{88}) and 
 (\ref{1-e1}), we have
\begin{eqnarray*}
|I_1[v_k]|&=&
 |\int_{ |y|\le R_k' }[(\Delta -\Delta_{g_k})(v_k-U)]
  (\nabla v_k\cdot y+\frac{n-2}2 v_k)\\
&\le &C \int_{ |y|\le R_k' } (|\bar b_i||\partial_i (v_k-U)|
+|\bar d_{ij}| |\partial_{ij}
(v_k-U)|)U(y)\\
&\le & C
 M_k^{ -\gamma_1 -\frac {6+\beta_2''}{n-2} }
\int_{ |y|\le R_k' }(1+|y|)^{3-n-\gamma_2}
\\
&&+CM_k^{ -\gamma_1-\frac 8{n-2} }
\int_{ |y|\le R_k' }
(1+|y|)^{ 4-n-\gamma_2}\\
&&+ C  \sum_{l=2}^3  M_k^{ -\gamma_1 -\frac {2l+\beta_l'''}{n-2} }
\int_{ |y|\le R_k' }(1+|y|)^{ -n+l-\gamma_2}.
\end{eqnarray*}
The following estimates are straight forward:
$$
|I_3[v_k]|=O(1) M_k^{ -\frac 8{n-2} }
(R_k')^{6-n},
\qquad
|I_4[v_k]|=O(1)   (R_k')^{-n},
$$
$$
|I_5[v_k]|
=O(1)  \int_{ |y|= R_k' }  (|\nabla U(y)|^2 R_k'
+U|\nabla U(y)|)=O(1)(R_k')^{2-n}.
$$

\vskip 5pt
\hfill $\Box$
\vskip 5pt

To prove  Theorem \ref{thm1} for $n\ge 6$ we need to
establish appropriate decay rates of the Riemannian
curvature tensor at the center of the conformal normal coordinate
system we use.

\begin{lem} For $3\le n\le 6$, let $\{v_k\}$ satisfy (\ref{ab1}), assume that for some 
 $\delta>0$,
\begin{equation}
v_k(y)\le C_1U(y)\quad \mbox{ for }\quad |y|\le 
\delta M_k^{ \frac 2{n-2} }.
\label{abcd1}
\end{equation}
Then for any
$\epsilon>0$ there exists
some constant $C$ independent of $k$ such that
for all $|y|\le \frac \delta 2$,
\begin{equation}
|\nabla^j(v_k-U)(y)|\le CM_k^{-2+\frac{2\epsilon}{n-2}}(1+|y|)^{-\epsilon-j},
\quad j=0,1,2.
\label{abcd2}
\end{equation}
\label{prop3n6}
\end{lem}

\begin{lem}
For $n\ge 7$, $0\le \bar a<n-6$, let $v_k$ satisfy (\ref{ab1}).
There exists
some small $\delta'>0$, depending only on $n$, such that
if
 $\{R_k'\}$ and $\{v_k\}$  satisfy
\begin{equation}
2\le R_k'
\le 2\delta' M_k^{ \frac 2 {n-2} },
\label{PP1}
\end{equation}
\begin{equation}
M_k^{ \frac 8{n-2} }
=O(1) (R_k')^{ 4+\bar a},
\label{PP5}
\end{equation}
and
\begin{equation}
v_k(y)\le CU(y), \qquad |y|\le R_k',
\label{Q2}
\end{equation}
then, for $j=0,1,2$, and for some $C'>0$ independent of $k$,
\begin{equation}
|\nabla^j(v_k-U)(y)|\le C' M_k^{ -\frac 8{n-2} }(1+|y|)^{6-n+\bar a-j},
\qquad  \forall\  |y|\le \frac 14  R_k'.
\label{ea10}
\end{equation}
\label{prop201}
\end{lem}
\begin{rem} Since we have established (\ref{eq4}) for $n\le 7$, we already know that
the hypothesis in Lemma \ref{prop201} is satisfied for $n=7$ with
$R_k'=\delta M_k^{ \frac 2{n-2} }$, $\bar a=0$, where
$\delta>0$ is some number independent of $k$.
\label{rem30}
\end{rem}

We prove Lemma  \ref{prop201} first and Lemma \ref{prop3n6} next.

\bigskip

\noindent{\bf Proof of Lemma  \ref{prop201}.}\
Let $w_k=v_k-U$, consider the equation for $w_k$:
\begin{equation}
\Delta_{g_k}w_k(y)-\bar cw_k(y)+n(n+2)\bar \xi^{\frac{4}{n-2}}w_k(y)=
\bar cU(y),
\label{54-1}
\end{equation}
where
$$
n(n+2)\bar \xi^{\frac{4}{n-2}}(y)=\left\{\begin{array}{ll}
\displaystyle{n(n-2)
\frac{v_k^{\frac{n+2}{n-2}}(y)-U^{\frac{n+2}{n-2}}(y)}{v_k(y)-U(y)}}&
\qquad \mbox{if}\quad v_k(y)\neq U(y), \\  \\
n(n+2)U(y)^{\frac{4}{n-2}}&\qquad \mbox{if}\quad  v_k(y)=U(y).
\end{array}
\right.
$$
By  (\ref{rough}),
\begin{equation}
|\bar cU(y)|\le  C
 M_k^{ -\frac
8 {n-2} } (1+|y|)^{4-n}.
\label{55-2}
\end{equation}
For $R_1$ sufficiently large,
  we claim that the operator
$\Delta_{g_k}-\bar c+n(n+2)\bar \xi^{\frac{4}{n-2}}$ satisfies maximum
principle over $R_1<|y|<R_k'$.  To see this, we estimate
the  $L^{\frac n2}$ norm of the coefficients of
$w_k$:
\begin{eqnarray*}
&&\int_{R_1\le |y|\le R_k' }|n(n+2)\bar \xi^{ \frac 4{n-2} }
-\bar c|^{ \frac n2}\\
&\le & C\int_{ R_1\le |y|\le \delta' M_k^{ -\frac 2{n-2} }  }
|M_k^{ -\frac 8{n-2} }|y|^2+|y|^{-4}|^{ \frac n2}
\le
 C[(\delta')^{2n}+ (R_1)^{-n}].
\end{eqnarray*}
So for $R_1$ sufficiently large and $\delta'$ sufficiently small, the maximum
principle holds for
$\Delta_{g_k}-\bar c+n(n+2)\bar \xi^{\frac{4}{n-2}}$ over
$R_1\le |y|\le R_k'$.

We compare $w_k$ with
$$
f(r):= C_{10} M_k^{ -\frac 8{n-2} }|y|^{6-n+\bar a}
$$
over $R_1\le |y|\le R_k'$
 where $C_{10}$ will be chosen momentarily.

Since
$\Delta (r^{6-n+\bar a})=
-(n-6-\bar a)(\bar a+4) r^{4-n+\bar a}$,
we have, using
 (\ref{rough}), (\ref{PP1}) and (\ref{Q2}), and taking $R_1$ larger
and $\delta'$ smaller if necessary,
\begin{eqnarray}
&&(\Delta+\bar b_i\partial_i+\bar d_{ij}\partial_{ij}-\bar c
+n(n+2)\bar \xi^{\frac{4}{n-2}})f(|y|)
\nonumber\\
&\le& -\frac 12 (n-6-\bar a)(\bar a+4) 
r^{-2}f(r),
\quad R_1\le |y|\le R_k'.
\label{PP3}
\end{eqnarray}
Making $C_{10}$ larger if necessary, we have, using (\ref{PP5})
and Corollary \ref{cor2-2new},
$$
|w_k(y)|\le f(y),\qquad |y|=R_1\ \mbox{or}\ |y|=R_k'.
$$
By the maximum principle, applied to the difference of
(\ref{54-1}) and (\ref{PP3}), and using (\ref{55-2}) and
 $\bar a\ge 0$, 
we have
$$
|w_k(y)|\le f(y),\qquad R_1<|y|\le R_k'.
$$
Estimate (\ref{ea10}) for $j=0$ follows from this.

Applying standard elliptic estimates to the equation of
$v_k$ and using (\ref{Q2}), we obtain
$$
|\nabla v_k(y)|\le (1+|y|)^{ 1-n},\
\qquad
|y|\le \frac 12 R_k'.
$$
With this and (\ref{59}), we obtain
\begin{equation}
|\nabla \bar \xi^{ \frac 4{n-2} }(y)|\le C(1+|y|)^{ -5},
\qquad
|y|\le \frac 12 R_k'.
\label{gg1}
\end{equation}
Let $\hat w_k(z)= w_k(y+\frac {|y|}{10}z)$,
$|z|\le 1$.
Applying standard elliptic estimates
to the equation satisfied by $\hat w_k$,
derived from (\ref{54-1}), and using (\ref{gg1}),
 we obtain
$$
|\nabla \hat w_k(0)|+|\nabla ^2 \hat w_k(0)|
\le C\bigg(\|\hat w_k\|_{ L^\infty(B_1) } +
 \| |y|^2 (\bar c U)(y+\frac { |y| }{10} \cdot)\|_{ C^1(B_1) }\bigg).
$$
Estimate  (\ref{ea10}) for $j=1,2$ follows from the above
by using (\ref{ea10}) for $j=0$.
Lemma  \ref{prop201} is proved.

\vskip 5pt
\hfill $\Box$
\vskip 5pt

\noindent{\bf Proof of Lemma \ref{prop3n6}:}
  This proof is essentially the same as that of Lemma \ref{prop201}.
We only need to change the comparison function to 
$$f(r):=C_{10}M_k^{-2+\frac{2\epsilon}{n-2}}r^{-2-\epsilon}$$ and
to observe that
$$|\bar cU(y)|\le CM_k^{-\frac{8}{n-2}}(1+|y|)^{4-n}
\le CM_k^{-2+\frac{2\epsilon}{n-2}}(1+|y|)^{-2-\epsilon}. \Box $$

\bigskip

\noindent{\bf Proof of (\ref{nn1}) for $n=6$.}\
By (\ref{pohov}) with $R_k' =\frac \delta 2 M_k^{ \frac 2{n-2} }$,
$$
I_2[v_k]=O(1)\left(|I_1[v_k]|+|I_3[v_k]|+|I_4[v_k]|+|I_5[v_k]|\right).
$$
Using (\ref{rough}), (\ref{88}), (\ref{abcd1}) and (\ref{abcd2}),
we obtain 
$$
|I_1[v_k]|+|I_3[v_k]|+|I_4[v_k]|+|I_5[v_k]|\le CM_k^{-2}.
$$
In view of the first line in (\ref{i2crude}), we have, for some
$c_6(n)>0$,
$$
I_2[v_k]=-12n c_6(n)\bar R^{ (2) }M_k^{-2}\log M_k
+O(M_k^{-2})=-c_6(n)|W|^2 M_k^{-2}\log M_k
+O(M_k^{-2}).
$$
Estimate  (\ref{nn1}) for $n=6$ follows from the above.

\medskip

\noindent{\bf Proof of (\ref{nn1}) for $n=7$.}\
By Lemma \ref{prop201} and Remark \ref{rem30},
(\ref{ea10}) holds for $n=7$, $R_k'=
\delta M_k^{ \frac 2{n-2} }$ and $\bar a=0$.
Thus, by Lemma \ref{cor2-6},  (\ref{2e-1})
holds with $\beta_2''=\beta_2'''=\beta_3'''=0$,
$R_k'=\delta M_k^{ \frac 2{n-2} }$,
$\gamma_1=\frac 8{n-2}$, $\gamma_2
=n-6$ and $\bar l=3$.  It follows that
$$
M_k^{ -\frac 8{n-2} }
\bar R^{(2)}\int_0^{ R_k'}
r^{1+n}U(r)^2dr
=O(M_k^{-2})
$$
which implies (\ref{nn1}) for $n=7$.

\begin{lem}
For $n\ge 8$, let  $v_k$
satisfy
(\ref{ab1}).  Then, for any $\epsilon>0$,
$$
|W|\le M_k^{ -\frac{2-\epsilon}{n-2} }.
$$
\label{lemweyl0}
\end{lem}

\noindent {\bf Proof of Lemma \ref{lemweyl0}.}\
By  Corollary \ref{cor2-5}, (\ref{Q2}) is satisfied 
with $R_k'=M_k^{ \frac{12-\epsilon}{  (n-2)^2 }  }$
for any small $\epsilon>0$.  With this $R_k'$,
(\ref{PP1}) is satisfied and  
(\ref{PP5}) is satisfied with $\bar a=\frac { 2(n-8) }{3} +O(1)\epsilon$. 
Thus by Lemma \ref{prop201},
(\ref{1-e1}) is satisfied with
$\gamma_1=\frac 8{n-2}$ and $\gamma_2=n-6-\bar a$.
Applying Lemma \ref{cor2-6} with the above data and $\bar l=3$,
we derive from (\ref{2e-1}) and (\ref{3-e0}) that
$$
|W|^2=-12n \bar R^{(2)}=O(1)
M_k^{ -\frac {4-\epsilon} {n-2} }.
$$
Lemma \ref{lemweyl0} is established.

\vskip 5pt
\hfill $\Box$
\vskip 5pt

The following properties of conformal normal
coordinates are established in \cite{HV}: If $W=0$, then
 $
  R_{abcd}=0$ and, for some constant
$c_1(n)>0$,
 $\bar R^{(4)}=
 -c_1(n)|R_{abcd,e}|^2$;
if $W=0$ and $\nabla W=0$, then
$R_{abcd,e}=0$.
Examining the proofs there, 
we arrive at
\begin{equation}
|R_{abcd}|=O(1)|W|,\quad |R_{abcd,e}|=O(|W|)+O(|\nabla_g W|),
\label{bb1}
\end{equation}
and
\begin{equation}
\bar R^{(4)}=-c_1(n)|R_{abcd,e}|^2 +O(|W|).
\label{second}
\end{equation}

It follows from  Lemma \ref{lemweyl0}, (\ref{bb1}) and (\ref{second}) that,
for any $\epsilon>0$, 
\begin{equation}
|R_{abcd}|=O(1)M_k^{ -\frac {2-\epsilon}{n-2} },
\label{first}
\end{equation}
and
\begin{equation}
\bar R^{(4)}\le CM_k^{ -\frac{2-\epsilon}{n-2} }.
\label{third}
\end{equation}

We know from Corollary \ref{prop5} that
\begin{equation}
\sigma_k\le CM_k^{ -\frac 8{n-2} }.
\label{M1}
\end{equation}
By (\ref{first}),
\begin{equation}
|\bar d_{ij}(y)|\le \left\{\begin{array}{ll}
CM_k^{ -\frac {6-\epsilon}{n-2} }|y|^{3-\frac \epsilon 2},&\quad |y|\ge 1,\\
CM_k^{ -\frac {6-\epsilon}{n-2} }|y|^2,&\quad |y|\le 1.
\end{array}
\right.
\label{M2}
\end{equation}

Following the proof of Proposition \ref{prop2new} while using
also (\ref{M1}) and (\ref{M2}), we have, instead of
(\ref{check}),
$$
E_\lambda(y)\le  \bar c(y) U^\lambda(y)-
(\frac \lambda {|y|})^{n+2}\bar c(y^\lambda)U(y^\lambda)
+O(1) M_k^{ -\frac{14-\epsilon}{n-2} }|y|^{1-n}.
$$
Instead of Corollary \ref{cor2-new2}, we now have
$$
E_\lambda(y)\le
 \sum_{l=2}^4 M_k^{ -\frac{4+2l}{n-2} }
H_{l,\lambda}(r)
\tilde R^{(l)}(\theta)
+ C_0
M_k^{ -\frac {14-\epsilon}{n-2} } r^{7-n-\frac \epsilon 2}.
$$
Take
$$
\tilde h_{1,\lambda}(y)=
 \sum_{l=2}^4\sum_{j=1}^l
\sum_{i=1}^{I_j}  \tilde h_{l,j,\lambda}^{(i)}(y),
$$
$$
\tilde h_{2,\lambda}(y)=
Q M_k^{ -\frac {14-\epsilon}{n-2} }
f_{n, \min\{n-7+\frac \epsilon 2, 2-\epsilon\} }(\frac r\lambda),
$$
and
$$
h_\lambda(y)=\tilde h_{1,\lambda}(y)+\tilde h_{2,\lambda}(y).
$$
Let
\begin{equation}
R_k=
\left\{
\begin{array}{rl}
k^{-\frac 14}M_k^{ \frac 2{n-2} },& n=8,\\
M_k^{ \frac{14-\sqrt{\epsilon}}{(n-2)^2} }, &n\ge 9.
\end{array}
\right.
\label{MM1}
\end{equation}
Then we can follow the proof of Lemma \ref{prop3} to show that
\begin{equation}
\min_{|y|=r}v_k(y)\le (1+\epsilon)U(r),\qquad\forall\  0<r\le
 R_k.
\label{eq35new}
\end{equation}
Indeed we only need to verify a few things.  First we still have
(\ref{Z1}). As before, we can show that  
$$
(\Delta_{g_k}-\bar c+n(n+2)\xi^{\frac{4}{n-2}}) 
\tilde h_{2,\lambda}(y)
\le -\frac Q4 \bar D_{k,\epsilon}(r),
$$
where 
$$
\bar D_{k,\epsilon}(r):=  M_k^{ -\frac {14-\epsilon}{n-2} }
r^{ -\min\{n-7+\frac \epsilon 2, 2-\epsilon\}}.
$$

We can verify, using the strengthened estimate
(\ref{M2}),  that
\begin{equation}
|\bar c||\tilde h_{1,\lambda}(y)|
+|\bar b_i \partial_i \tilde h_{1,\lambda}(y)|
+|\bar d_{ij}\partial_{ij}  \tilde h_{1,\lambda}(y)|
\le C \bar D_{k,\epsilon}(r).
\label{MM7}
\end{equation}

Recall that  
\begin{equation}
|\nabla^j(v_k-U)(y)|\le C' M_k^{ -\frac 8{n-2} }|y|^{6-n+\bar a-j},
\qquad  \forall\  |y|\le  \frac 14  M_k^{ \frac{12-\epsilon}{(n-2)^2}  },
\label{MM5}
\end{equation}
where $\bar a=\frac{2(n-8)}3 +\sqrt{\epsilon}$.
With (\ref{MM5}) we have, instead of
(\ref{est1}),
\begin{equation}
|n(n+2)\xi^{\frac{4}{n-2}}(y)-V_{\lambda}(|y|)| |\tilde h_{1,\lambda}(y)|
\le
CM_k^{-\frac{8}{n-2}}
|y|^{\bar a}, \qquad
\lambda\le |y|\le  \frac 14 M_k^{ \frac{12-\epsilon}{(n-2)^2}  },
\label{MM6}
\end{equation}
which can be shown by following
the  arguments in the proof of Lemma \ref{lem8-1} together
with the improved bounds 
$$
|a(y)|+|b(y)|\le M_k^{ -\frac 8{n-2} } |y|^{6-n+\bar a}
$$
given by (\ref{MM5}).

With (\ref{est2}) and the improved estimate (\ref{MM6}),
we can show that
\begin{eqnarray*}
&&
(\Delta+n(n+2)\xi^{\frac 4{n-2} })\tilde h_{1,\lambda}
\\
&\le & (\Delta+V_\lambda)\tilde h_{1,\lambda}
+|n(n+2)\xi^{\frac 4{n-2} }-V_\lambda||\tilde h_{1,\lambda}|
\\
&\le &
- \sum_{l=2}^4  M_k^{ -\frac{4+2l}{n-2} }
H_{l,\lambda}(r) \tilde R^{(l)}(\theta)
+ O(1) \bar  D_{k,  \epsilon}(r),
\qquad \mbox{in}\  \hat O_\lambda.
\end{eqnarray*}
Thus, in view of 
(\ref{MM7}), 
\begin{eqnarray}
&&
(\Delta_{g_k}-\bar c+n(n+2)\xi^{\frac 4{n-2} })\tilde h_{1,\lambda}
\nonumber\\
&\le& - \sum_{l=2}^4 M_k^{ -\frac{4+2l}{n-2} }
H_{l,\lambda}(r) \tilde R^{(l)}(\theta)
+C \bar  D_{k,  \epsilon}(r),
\qquad
\mbox{in}\  \hat O_\lambda.
\nonumber
\end{eqnarray}
Fixing a large $Q$, we obtain
\begin{equation}
(\Delta_{g_k}-\bar c+n(n+2)\xi^{\frac{4}{n-2}})h_{\lambda}+E_{\lambda}\le 0
\quad \mbox{in}\quad \hat O_{\lambda}.
\label{eq19new}
\end{equation}
With (\ref{eq19new}), estimate (\ref{eq35new}) with 
$R_k$ given by (\ref{MM1}) is established, as in the proof of  
Lemma \ref{prop3}.

Once (\ref{eq35new})  is established, the proof of Lemma \ref{prop4}
yields
\begin{equation}
v_k(y)\le CU(y),\qquad |y|\le  R_k
\label{MM8}
\end{equation}
for the $R_k$ given by  (\ref{MM1}).

\noindent{\bf Proof of
(\ref{eq4}) for $n=8$.}\ Since estimate (\ref{MM8}) holds
for $R_k=k^{ -\frac 12}M_k^{ \frac 2{n-2} }$, which
violates (\ref{200}), estimate (\ref{eq4}) in dimension
$n=8$  is established.

Now we turn to $n\ge 9$.  Since
(\ref{MM8}) holds with $R_k=M_k^{ \frac{14-\epsilon}{(n-2)^2 }  }$
for any $\epsilon>0$,  
hypotheses (\ref{PP1}), (\ref{PP5}) and
(\ref{Q2}) hold with $R_k' =M_k^{ \frac{14-\epsilon}{(n-2)^2 }  }$
and $\bar a=\frac { 4(n-9) }7+\sqrt{\epsilon}$.  Thus, by Lemma
\ref{prop201},
\begin{equation}
|\nabla^j(v_k-U)(y)|\le C' M_k^{ -\frac 8{n-2} }|y|^{6-n+\bar a-j},
\qquad  \forall\  |y|\le  M_k^{ \frac{14-\epsilon}{(n-2)^2 }  },
\ j=0,1,2.
\label{abc1}
\end{equation}
With this, we can apply Lemma \ref{cor2-6} with
$\bar l=5$, $\beta_2''=\beta_3'''=
0$, $\beta_2'''=2-\epsilon$,
 $R_k'=  M_k^{ \frac{14-\epsilon}{(n-2)^2 }  }$,
$\gamma_1= \frac 8{n-2}$ and 
$\gamma_2=n-6-\bar a$ to deduce,
in view of (\ref{second}),  from (\ref{2e-1}) and
(\ref{3-e0}) that
$$
M_k^{ -\frac 8{n-2} }|W|^2+
M_k^{ -\frac{12}{n-2} }\big[ |R_{abcd,e}|^2+
O(|W|)\big]
=O(M_k^{ -\frac{14-\epsilon}{n-2} }).
$$
Consequently, we have, using (\ref{bb1}), that for any $\epsilon>0$,
\begin{equation}
|W|=O(1)M_k^{ -\frac {3-\epsilon}{n-2} },
\qquad |R_{abcd,e}|=O(1) M_k^{ -\frac {1-\epsilon}{n-2} }. 
\label{NN1}
\end{equation}

\bigskip

Now we consider $n\ge 9$.
By (\ref{NN1}),
$$
|b_i(x)|\le CM_k^{ -\frac{\beta_2''}{n-2} }|x|^2
+C|x|^3,\ 
|d_{ij}(x)|\le C\sum_{l=2}^3
M_k^{ -\frac{\beta_l'''}{n-2} }|x|^l+
C|x|^4,
$$
or, equivalently,
\begin{equation}
\left\{
\begin{array}{l}
|\bar b_i(y)|\le 
 CM_k^{ -\frac{6+\beta_2''}{n-2} }|y|^2+
CM_k^{ -\frac 8{n-2} }|y|^3,
\\ 
|\bar d_{ij}(y)|\le 
 C\sum_{l=2}^3
M_k^{ -\frac{2l+\beta_l'''}{n-2} }|y|^l+
CM_k^{ -\frac 8{n-2} }|y|^4,
\end{array}
\right.
\label{AB2}
\end{equation}
where
$$
\beta_2''=\beta_3'''=1-\epsilon, \ \beta_2'''=3-\epsilon,
\ \epsilon>0.
$$
With (\ref{AB2}),
we follow the proof of
Proposition \ref{prop2new} to obtain, using also (\ref{M1}),
that
$$
E_\lambda(y)= \bar c(y) U^\lambda(y)-
(\frac \lambda {|y|})^{n+2}\bar c(y^\lambda)U(y^\lambda)
+O(1)M_k^{ -\frac{15-\epsilon}{n-2} }|y|^{2-n}.
$$
We also know that for $n\ge 9$,
$$
\bar R^{(4)}\le C M_k^{ -\frac{3-\epsilon}{n-2} }.
$$
Thus we  have, using also $\bar R^{(2)}\le 0$ and
(\ref{RR2}) with $\bar l=5$,
\begin{equation}
E_\lambda(y)\le
 \sum_{l=2}^5 M_k^{ -\frac{4+2l}{n-2} }
H_{l,\lambda}(r)
\tilde R^{(l)}(\theta)
+ CM_k^{ -\frac{15-\epsilon}{n-2} }
|y|^{-\frac 32 -\frac \epsilon2}. 
\label{BB1}
\end{equation}

Now we take
$$
\tilde h_{1,\lambda}(y)=
 \sum_{l=2}^5\sum_{j=1}^l
\sum_{i=1}^{I_j}  \tilde h_{l,j,\lambda}^{(i)}(y),
$$
$$
\tilde h_{2,\lambda}(y)=
Q M_k^{ -\frac {15-\epsilon}{n-2} }
f_{n, \frac 32 +\epsilon}
(\frac r\lambda),
$$
and
$$
h_\lambda(y)=\tilde h_{1,\lambda}(y)+\tilde h_{2,\lambda}(y).
$$
Let
\begin{equation}
R_k=
k^{-\frac 14}M_k^{ \frac 2{n-2} }.
\label{MM1new}
\end{equation}
Then we can follow the proof of Lemma \ref{prop3} to show,
with the above $\{R_k\}$,  that
\begin{equation}
\min_{|y|=r}v_k(y)\le (1+\epsilon)U(r),\qquad\forall\  0<r\le
 R_k.
\label{eq35newnew}
\end{equation}
Indeed  we only need to verify a few things.  First we still have
(\ref{Z1}).  As before, we can show that
$$
(\Delta_{g_k}-\bar c+n(n+2)\xi^{\frac{4}{n-2}})
\tilde h_{2,\lambda}(y)
\le -\frac Q4 \bar D_{k,\epsilon}(r),
$$
where 
$$
\bar D_{k,\epsilon}(r):=  M_k^{ -\frac {15-\epsilon}{n-2} }
r^{ -\frac 32 - \frac \epsilon 2}.
$$

We now have the strengthened estimates:
$$
|\bar b_i(y)|\le C M_k^{ -\frac {7-\epsilon}{n-2} }
|y|^{ \frac {5-\epsilon}2 },\quad
|\bar d_{ij}(y)|\le C  M_k^{ -\frac {7-\epsilon}{n-2} }
|y|^3.
$$
With this we can show that
\begin{equation}
|\bar c||\tilde h_{1,\lambda}(y)|
+|\bar b_i \partial_i \tilde h_{1,\lambda}(y)|
+|\bar d_{ij}\partial_{ij}  \tilde h_{1,\lambda}(y)|
\le C \bar D_{k,\epsilon}(r).
\label{MM7new}
\end{equation}

Because of  (\ref{abc1}), we have, instead of
(\ref{est1}),
\begin{equation}
|n(n+2)\xi^{\frac{4}{n-2}}(y)-V_{\lambda}(|y|)| 
\le
CM_k^{-\frac{8}{n-2}}
|y|^{\bar a }, \qquad
\lambda\le |y|\le  \frac 14 M_k^{ \frac{14-\epsilon}{(n-2)^2}  },
\label{MM6new}
\end{equation}
which can be shown by following
the  arguments in the proof of Lemma \ref{lem8-1} together
with the improved bounds
$$
|a(y)|+|b(y)|\le M_k^{ -\frac 8{n-2} } |y|^{6-n+\bar a}
$$
given by  (\ref{abc1}).  Recall that $\bar a
=\frac { 4(n-9)+\sqrt{\epsilon}}  { 7}$.

With (\ref{est2}) and the improved estimate (\ref{MM6new}),
we can show that
\begin{eqnarray*}
&&
(\Delta+n(n+2)\xi^{\frac 4{n-2} })\tilde h_{1,\lambda}
\\
&\le & (\Delta+V_\lambda)\tilde h_{1,\lambda}
+|n(n+2)\xi^{\frac 4{n-2} }-V_\lambda||\tilde h_{1,\lambda}|
\\
&\le &
- \sum_{l=2}^5 M_k^{ -\frac{4+2l}{n-2} }
H_{l,\lambda}(r) \tilde R^{(l)}(\theta)
+ O(1) \bar  D_{k,  \epsilon}(r),
\qquad \mbox{in}\  \hat O_\lambda.
\end{eqnarray*}
Thus, in view of
(\ref{MM7new}),
\begin{eqnarray}
&&
(\Delta_{g_k}-\bar c+n(n+2)\xi^{\frac 4{n-2} })\tilde h_{1,\lambda}
\nonumber\\
&\le& - \sum_{l=2}^5 M_k^{ -\frac{4+2l}{n-2} }
H_{l,\lambda}(r) \tilde R^{(l)}(\theta)
+C \bar  D_{k,  \epsilon}(r),
\qquad
\mbox{in}\  \hat O_\lambda.
\nonumber
\end{eqnarray}
Fixing a large $Q$, we obtain
\begin{equation}
(\Delta_{g_k}-\bar c+n(n+2)\xi^{\frac{4}{n-2}})h_{\lambda}+E_{\lambda}\le 0
\quad \mbox{in}\quad \hat O_{\lambda}.
\label{eq19newnew}
\end{equation}
With (\ref{eq19newnew}), estimate (\ref{eq35newnew}) with
$R_k$ given by (\ref{MM1new}) is established, as in the proof of
Lemma \ref{prop3}.

Once (\ref{eq35newnew})  is established, the proof of Lemma \ref{prop4}
yields
\begin{equation}
v_k(y)\le CU(y),\qquad |y|\le  R_k
\label{MM8new}
\end{equation}
for the $R_k$ given by  (\ref{MM1new}).

\noindent{\bf Proof of
(\ref{eq4}) for $n=9$.}\ Since estimate (\ref{MM8new}) holds
for $R_k=k^{ -\frac 12}M_k^{ \frac 2{n-2} }$, which
violates (\ref{200}), estimate (\ref{eq4}) in dimension
$n=9$  is established. $\Box$

\bigskip

Now we make use of the Pohozaev identity (\ref{pohov}) 
(with appropriate $R_k'$) to prove (\ref{nn1}), (\ref{nn2}) and (\ref{mm1}).
Since we have established (\ref{eq4}), the assumption in Lemma
\ref{prop3n6} is satisfied for some $\delta>0$ and the assumption of 
Lemma \ref{prop201} is satisfied for $7\le n\le 9$ for $\bar a=0$ and
for some $\delta'>0$ Thus, for some $\delta>0$, we have , for some
$\epsilon>0$, any $|y|\le \delta M_k^{\frac{2}{n-2}}$,
\begin{equation}
|\nabla^j(v_k-U)(y)|\le C(\delta,\epsilon)M_k^{-2+\frac{2\epsilon}{n-2}}
(1+|y|)^{-\epsilon-j},j=0,1,2,\quad 3\le n\le 6,
\label{eq170}
\end{equation}
\begin{equation}
|\nabla^j(v_k-U)(y)|\le C(\delta,\epsilon)M_k^{-\frac{8}{n-2}}(1+|y|)^{6-n-j},
j=0,1,2,\quad 7\le n\le 9.
\label{eq173}
\end{equation}

Taking
$R_k'=\sigma M_k^{\frac{2}{n-2}}$ in (\ref{pohov}), $0<\sigma<\delta$, using
(\ref{eq170}) and (\ref{rough}), we have, for some $C_1(n)>0$,
\begin{equation}
|M_k^2I_2[v_k]-C_1(n)|W|^2M_k^{\frac{2(n-6)}{n-2}}\int_{|y|\le \sigma M_k^{\frac{2}{n-2}}}
|y|^2U^2(y)|\le
\left\{\begin{array}{ll}
C\sigma^2,&\quad n=6,\\
C\sigma, &\quad n=7.
\end{array}
\right.
\label{A2-1new}
\end{equation}

\begin{equation}
|M_k^2I_1[v_k]|\le
\left\{\begin{array}{ll}
C\sigma^{2-\epsilon},&\quad n=6,\\
C\sigma, &\quad n=7.
\end{array}
\right.
\label{A3-1new}
\end{equation}

\begin{equation}
|M_k^2I_3[v_k]|\le C,\quad 6\le n\le 9.
\label{A3-2new}
\end{equation}

\begin{equation}
\limsup_{k\to \infty} M_k^2|I_4[v_k]|=0,\quad 6\le n\le 9.
\label{A3-3new}
\end{equation}

\begin{equation}
\limsup_{k\to \infty} M_k^2|I_5[v_k]|=0, \quad 6\le n\le 9.
\label{A3-4new}
\end{equation}

Estimate (\ref{nn1}) for $n=6,7$ follows from (\ref{pohov}), (\ref{A2-1new}),
(\ref{A3-1new}),(\ref{A3-2new}),(\ref{A3-3new}) and (\ref{A3-4new}).
By Lemma \ref{lemweyl0}, (\ref{bb1}), (\ref{NN1}), assumption (\ref{TT1}) holds
for
$$\beta_2''=\beta_3'''=\left\{\begin{array}{ll}0,&\quad n=8,\\
1-\epsilon,&\quad n=9,
\end{array}
\right.
$$
$$\beta_2'''=\left\{\begin{array}{ll}2-\epsilon,&\quad n=8,\\
3-\epsilon,&\quad n=9,
\end{array}
\right.
$$
where $\epsilon>0$ is any number. Let $R_k'=\sigma M_k^{\frac{2}{n-2}}$,
$0<\sigma<\delta$, $\gamma_1=\frac{8}{n-2}$, $\gamma_2=n-6$, $\bar l=5$,
we deduce from (\ref{2e-1}), in view of (\ref{second}), that
$$|W|^2M_k^{\frac{2(n-6)}{n-2}}+|R_{abcd,e}|^2
M_k^{\frac{2(n-8)}{n-2}}\int_0^{\sigma M_k^{\frac{2}{n-2}}}r^{3+n}U^2(r)dr
=O(1)$$
for $n=8,9$. 
Estimate (\ref{nn1}) and (\ref{nn2})  for $n=8,9$ follows from the above.
\bigskip

\noindent{\bf Proof of (\ref{mm1}).}\
For $n\ge 10$, we claim that we have
\begin{eqnarray}
&&v_k(y)\le CU(y),\quad |y|\le R_k'=M_k^{\frac{16-\epsilon}{(n-2)^2}}.
\label{mm11}\\
&&\beta_2''=\beta_3'''=2-\epsilon,\quad \beta_2'''=4-\epsilon.
\label{oct7e1}
\end{eqnarray}
Given this and (\ref{ea10}) with $\bar a=
\frac{ n-10 +\sqrt{\epsilon} }{2}$, an application of
 Lemma \ref{cor2-6} yields
$$c_1(n)|W(x_k)|^2M_k^{\frac{2(n-6)}{n-2}}+c_2(n)
|\nabla W(x_k)|^2M_k^{\frac{2(n-8)}{n-2}}=O(M_k^2(R_k')^{2-n}).$$
Then estimate (\ref{mm1}) follows from the above.

To show (\ref{mm11}), we need to improve the decay rate of $|W(x_k)|$
and $|\nabla W(x_k)|$ by iteration. Right now for the rate of
$\beta_2''$, $\beta_2'''$ and $\beta_3'''$ we have (\ref{AB2}),
for $E_{\lambda}$ we have (\ref{BB1}).
By exactly the same way of constructing auxiliary
functions we have, for $n\ge 10$, that
\begin{equation}
v_k(y)\le CU(y),\quad |y|\le R_k':=M_k^{\frac{15-5\epsilon}{(n-2)^2}}.
\label{mm12}
\end{equation}
Since $\epsilon$ indicates an arbitrarily small positive number, we just
replace $5\epsilon$ by $\epsilon$. Then we apply Lemma \ref{cor2-6}
with this $R_k'$ to get
$$c_1(n)|W|^2M_k^{-\frac{8}{n-2}}+c_2(n)|\nabla W|^2M_k^{-\frac{12}{n-2}}
=O((R_k')^{2-n})=O(1)M_k^{-\frac{15-\epsilon}{n-2}}.$$
So after this step we have
$$ |W|=O(M_k^{-\frac{3.5-\epsilon}{n-2}}),\quad  |\nabla W|
=O(M_k^{-\frac{1.5-\epsilon}{n-2}}),\quad \mbox{for any }\epsilon>0.$$

By this stronger decay rate we can show (\ref{mm12}) for
$R_k'=M_k^{\frac{15.5-\epsilon}{(n-2)^2}}$. Then Pohozaev identity gives
$|W|=O(M_k^{-\frac{3.75-\epsilon}{n-2}})$ and
$|\nabla W|=O(M_k^{-\frac{1.75-\epsilon}{n-2}})$. The corresponding rates
for $\beta_2'',\beta_2''',\beta_3'''$ improve too. These new rates lead
to (\ref{mm11}) except that $R_k'=M_k^{\frac{15.75-\epsilon}{(n-2)^2}}$.  
Now we see that for any
fixed $\epsilon'>0$, after doing this iteration finite times we can 
derive (\ref{mm11}) for $R_k'=M_k^{\frac{16-\epsilon'}{(n-2)^2}}$.
(\ref{mm1}) is established. So is Theorem \ref{thm1}.  $\Box$

\section{Proof of Theorem \ref{thm0} for  $p=\frac{n+2}{n-2}$}
In this section we prove  Theorem \ref{thm0} for $p=\frac{n+2}{n-2}$ 
by contradiction argument. By Remark \ref{rem10},
we only need to consider $3\le n\le 9$.
Suppose the contrary, for each $k=1,2,\cdots,$
$C=k$ does not work for some solution $u_k$ of
(\ref{Y0}) satisfying the hypotheses of the theorem.
Clearly $u_k(\bar P)\to \infty$.
A point $Q$ in $\Omega$ is called a
blow up point of $\{u_k\}$ if 
$\{u_k\}$  is not locally bounded in some fixed neighborhood
of $Q$.  Let ${\cal S}^*$ denote the set of all blow up
points of $\{u_k\}$ in $\Omega$.

\begin{lem}  There exists some positive constant
$\delta'$, depending only on
$(M,g)$, $\Omega$, $dist_g(\bar P, \partial \Omega)$,
$\bar b$ and $\epsilon$ such that
$$
dist_h(Q, Q')\ge \delta',\qquad
\forall\ Q, Q'\in 
{\cal S}^*\cap \Omega_\epsilon, Q\ne Q'.
$$
\label{lem11-1}
\end{lem}

\noindent{\bf Proof of Lemma \ref{lem11-1}.}\  
$\forall\ Q\in ( {\cal S}^*\cap \Omega_\epsilon)\setminus\{\bar P\}$, there
exists a subsequence of $\{u_k\}$, still denoted as
$\{u_k\}$, and a sequence of points
$P_k'\to Q$ such that
$u_k(P_k')\to \infty$.  Let
$\delta_k=u_k(P_k')^{ -\frac 1{n-2}}$, then
$\delta_k^{  \frac {n-2}2 } u_k(P_k')\to \infty$. 

Consider 
$$
\tilde u_k(P)=(\delta_k-dist_g(P, P_k'))^{ \frac {n-2}2  }
u_k(P), \qquad
P\in B(P_k', \delta_k),
$$
and let $P_k''$ be a maximum point of $\tilde u_k$ in the
closure of $B(P_k', \delta_k)$.
Let $r_k=\frac 12(\delta_k-dist_g(P_k'', P_k'))
\in (0, \frac{\delta_k}2).
$
Then
\begin{eqnarray}
\gamma_k^{ \frac 2{n-2} }&:=&
(r_k)^{ \frac {n-2}2 }u_k(P_k'')
=2^{ \frac {2-n}2 }
\tilde u_k(P_k'')\ge 2^{  \frac {2-n}2 }\tilde u_k(P_k')
\nonumber\\
&=&  2^{  \frac {2-n}2 } (\delta_k)^{ \frac {n-2}2 } u_k(P_k')\to\infty,
\label{E4-1}
\end{eqnarray}
and
$$
(2r_k)^{ \frac {n-2}2 } u_k(P_k'')
= \tilde u_k(P_k'')\ge \tilde u_k(P)
\ge (r_k)^{ \frac {n-2}2 } u_k(P),\quad
\forall\ P\in B(P_k'', r_k),
$$
i.e.
\begin{equation}
\sup_{ B(P_k'', r_k)  }u_k
\le 2^{ \frac {n-2}2 } u_k(P_k'').
\label{E4-2}
\end{equation}

Let $\{x^1, \cdots, x^n\}$ be some geodesic normal coordinates centered at $P_k''$, so $x=0$ corresponds to
$P_k''$.  Consider
$$
w_k(y)=\frac 1{ u_k(0)}
u_k\left( \frac y{ u_k(0)^{ \frac 2{n-2} }  }\right),
\qquad |y|<\Gamma_k:=
\frac 12 dist_g(Q, \partial \Omega)
u_k(0)^{ \frac 2{n-2} }.
$$
Then, with $g_k$ denoting the rescaled metric,
$$
-L_{g_k}w_k=n(n-2) w_k^{ \frac{n+2}{n-2}},\qquad
|y|<\Gamma_k,
$$
$$
1=w_k(0)\ge 2^{ \frac {2-n}2 }w_k(y),
\qquad |y|<\gamma_k,
$$
$$
\Gamma_k\ge \gamma_k\to\infty.
$$
As usual, $w_k\to w$ in $C^2_{loc}(\Bbb R^n)$ with
$$
w(y)=\left( \frac {\bar \lambda}{ 1+\bar \lambda^2|y-\bar y|^2}\right)^{ \frac {n-2}2},
$$
and $0<\bar \epsilon\le \bar\lambda\le 1/\bar\epsilon$, 
$|\bar y|\le 1/\bar\epsilon$ for some constant $\bar \epsilon$ depending only
on $n$.
It follows that $\nabla w_k(y)=0$ for some
$y_k=\bar y+\circ(1)$.
Let $P_k'''=u_k(P_k'')^{ -\frac 2{n-2} }y_k$, then
\begin{equation}
u_k(P_k'')\le u_k(P_k''')\le (2^{ \frac {n-2}2}+\circ(1))u_k(P_k''),
\qquad
\nabla u_k(P_k''')=0.
\label{abcde1}
\end{equation}
By the convergence of $w_k$ to $w$, we have, for some
$\hat \epsilon>0$ depending only on $n$,
\begin{equation}
u_k(P)\ge \hat \epsilon u_k(P_k'''),
\qquad \forall\ dist_g(P, P_k''')
=u_k(P_k''')^{ -\frac 2{n-2} }.
\label{E6-1}
\end{equation}
Let $G$ denote the Green$'$s function of $-L_g$ on
$\Omega$ with zero Dirichlet boundary condition.  Then
$G$ is positive in $\Omega$ and
\begin{equation}
\lim_{ k\to\infty}\max_{
dist_g(P_k''', P)=u_k(P_k''')^{ -\frac 2{n-2} }
}
|G(P_k''', P)dist_g(P_k''',P)^{n-2}
-\frac 1{  (n-2) \omega_n }|=0,
\label{E7-1}
\end{equation}
where $\omega_n$ denotes the volume of the standard
$(n-1)-$sphere.

By (\ref{E6-1}) and (\ref{E7-1}),
$$
u_k(P)\ge [(n-2)\omega_n \hat \epsilon+\circ(1)]
u_k(P_k''')^{-1} G(P_k''', P),
\quad\forall\ dist_g(P_k''',P)=u_k(P_k''')^{-\frac 2{n-2} }.
$$
Since $L_gu_k\le0$, we have,
using the maximum principle,
\begin{equation}
u_k(P)\ge  [(n-2)\omega_n \hat \epsilon+\circ(1)]
u_k(P_k''')^{-1} G(P_k''', P)
\label{E8-1}
\end{equation}
for all $P\in \Omega$ satisfying $\mbox{dist}_g(P_k''',P)
\ge u_k(P_k''')^{-\frac 2{n-2} }$.
By Theorem \ref{thm1}, there exist
some universal constants $C_0, \delta>0$, i.e., they
depend only on $(M,g)$,
$\Omega, dist_g(\bar P, \partial \Omega)$,
$\bar b$ and $\epsilon$, such that
\begin{equation}
u_k(\bar P)u_k(P)dist_g(\bar P, P)^{n-2}
\le C_0,\qquad \forall\ 0<dist_g(\bar P,P)\le \delta.
\label{E9-1}
\end{equation}
Taking a $P$ in $B(\bar P, \delta)\setminus B(\bar P, \frac \delta 2)$ 
satisfying $dist_g(P_k''',P)\ge \frac \delta 9$, we obtain, using (\ref{E9-1}) and
(\ref{E8-1}),
$$
u_k(\bar P)\le b' u_k(P_k'''),
$$
where $b'\ge 1$ is some universal constant.  Thus,
by (\ref{Y1}), 
$$
\sup_\Omega u_k
\le \bar b b' u_k(P_k''').
$$
Now, applying again Theorem \ref{thm1}, we have, for some
universal constants $C_0', \delta'>0$,
$$
u_k(P_k''')u_k(P) dist_g(P_k''',P)^{n-2}
\le C_0',\qquad
\forall\ 0<dist_g(P_k'',P)\le 2\delta'.
$$
It follows, using the fact $P_k'''\to Q$, that
$$
{\cal S}^*\cap B(Q, \delta')=\{Q\}.
$$
Lemma \ref{lem11-1} is established.

\vskip 5pt
\hfill $\Box$
\vskip 5pt

Let ${\cal S}^*\cap \Omega_\epsilon=\{P_1, \cdots, P_m\}$.
Because of Lemma \ref{lem11-1}, we may choose local maximum
points $\{P_1^{(k)}, \cdots, P_m^{(k)}\}$ of $u_k$
such that $u_k(P_i^{(k)})\to \infty$ and
$P_i^{(k)}\to P_i$.   By the arguments
in the proof of Lemma \ref{lem11-1},
\begin{equation}
0<\liminf_{k\to\infty}
\frac{  u_k(P_i^{(k)})  }
{u_k(\bar P)  }
\le \limsup _{k\to\infty}
\frac{  u_k(P_i^{(k)})  }
{u_k(\bar P)  }<\infty,
\quad 1\le i\le m.
\label{E11-1}
\end{equation}
By Theorem \ref{thm1}, there exist
some positive constants $C_0>0$ and
$0<\delta<\frac \epsilon 2$ such that
\begin{equation}
u_k(P_i^{(k)}) u_k(P)
dist_g(P_i^{(k)},P)\le C_0,\quad
\forall\ dist_g(P_i^{(k)},P)\le \delta,\ 1\le i\le m.
\label{E11-2}
\end{equation}

Since $\{u_k\}$
is bounded in $\Omega_{2\epsilon}\setminus
\cup_{ i=1}^m B(P_i^{(k)}, \frac \delta 2)$, we have, by the Harnack inequality and (\ref{E11-2}),
\begin{equation}
u_k(P)\le Cu_k(\bar P)^{-1},
\qquad \forall\ P\in
\Omega_{2\epsilon}\setminus \cup_{i=1}^m B(P_i^{(k)}, \delta).
\label{E12-1}
\end{equation}

Now applying Theorem \ref{thm1} together with some
standard elliptic estimates, we obtain
(\ref{ee7}), (\ref{2-5}), (\ref{nn1new}),
(\ref{nn2new}) and (\ref{2-6}), with $u$ and $\{P_i\}$ replaced
respectively by $u_k$
and $\{P_i^{(k)}\}$,
for some constant $C$ independent of $k$.
Estimate (\ref{abc}), still with 
  $u$ and $\{P_i\}$ replaced 
  respectively by $u_k$
  and $\{P_i^{(k)}\}$, for  some constant $C$ independent of $k$, 
 is simply a rewritten of the estimates on the rescaled 
 $v_k$ obtained in Section 2.
 Theorem \ref{thm0}
is established.

\vskip 5pt
\hfill $\Box$
\vskip 5pt

\section{Proof of Theorem \ref{thm3} for $p=\frac{n+2}{n-2}$}

In this section we establish Theorem \ref{thm3} for $p=\frac{n+2}{n-2}$.
The case $1<1+\epsilon\le p\le \frac{n+2}{n-2}$ will be discussed
in Section 5.
Suppose the contrary of
(\ref{compact}), 
let $\{u_k\}$ be a sequence of solutions of  (\ref{Y0})
with $p=\frac {n+2}{n-2}$ satisfying
\begin{equation}
u_k(P^{(k)})=\max_{M} u_k\to \infty,
\label{PK}
\end{equation}
with $P^{(k)}\to \bar P\in M$.

By Theorem \ref{thm0}, there exist local maximum points 
$\{P_1^{(k)}, \cdots, P_m^{(k)}\}$ of $u_k$,
such that (\ref{ee7}), (\ref{2-5}),     (\ref{nn1new}), (\ref{nn2new}),
 (\ref{2-6}), 
and (\ref{abc}) hold with 
$P_i$ replaced by $P_i^{(k)}$.
We may assume without loss of generality that $P_1^{(k)}=P^{(k)}$ and 
$P_i^{(k)}\to \bar P_i$ as $k\to\infty$.
As explained before we may
work in conformal
normal coordinates centered at $P_1^{(k)}=P^{(k)}$, and 
we rescale $u_k$ to $v_k$ as in (\ref{vk}).
Multiplying the equation of $u_k$ by $u_k(0)$, we obtain 
\begin{equation}
u_k(0)u_k\to G:=\sum_{i=1}^m a_i G(\cdot, \bar P_i),\qquad
\mbox{in}\ C_{\mbox{loc}}^2(M\setminus\{\bar P_1, \cdots, \bar P_m\}),
\label{nn5}
\end{equation}
where $a_i>0$.

By (\ref{nn1new}) and (\ref{nn2new})
\begin{equation}
W_g(\bar P_1)=0, \qquad \mbox{if}\ n\ge 6,
\label{nn7}
\end{equation}
and
$$\nabla _g W_g(\bar P_1)=0, \qquad \mbox{if}\ n\ge 8.$$

We take a small ball $B_{\sigma}:=B(\bar P_1,\sigma)$.
The Pohozaev identity is computed over $\bar B_{\sigma}$.
The following lemma will be combined with the positive mass theorem
to get a contradiction.

\begin{lem}
$$\lim_{\sigma\to 0+}B(\sigma,G,\nabla G)\ge 0,\quad 3\le n\le 9.$$
\label{Blem}
\end{lem}

\noindent{\bf Proof of Lemma \ref{Blem}:}

Since  $u_k(0)u_k\to G$ in $C_{\mbox{loc}}^2(B_{\sigma}\setminus \{0\})$,
$M_k^2I_5[v_k]\to B(\sigma,G,\nabla G)$ over $\partial B_{\sigma}$.
Thus by (\ref{pohov}) we only need to show
\begin{equation}
\lim_{\sigma\to 0+}\liminf_{k\to \infty}
M_k^2(I_1[v_k]+I_2[v_k]+I_3[v_k]+I_4[v_k])\ge 0. 
\label{39-1new}
\end{equation}

By (\ref{nn5}), 
\begin{equation}
\lim_{k\to \infty}M_k^2I_4[v_k]=-\frac{(n-2)^2}2\sigma \lim_{k\to \infty}
u_k(0)^2\int_{|x|=\sigma}u_k^{\frac{2n}{n-2}}=0,\quad 3\le n\le 9.
\label{C1-1}
\end{equation}

By (\ref{rough}), (\ref{rrr}), (\ref{eq170}) and (\ref{eq173})
\begin{equation}
M_k^2|I_1[v_k]|\le CM_k^2\int_{|y|\le \sigma M_k^{\frac{2}{n-2}}}
|(\bar b_i\partial_i+\bar d_{ij}\partial_{ij})(v_k-U)|U\le C\sigma,
\quad 3\le n\le 7.
\label{C1-2}
\end{equation}

By (\ref{SS5})
\begin{eqnarray}
&&M_k^2I_2[v_k] \nonumber \\
&=& -\frac{c(n)}{2}M_k^{2-\frac{4}{n-2}}
\int_{|y|\le \sigma M_k^{\frac{2}{n-2}}}\{2R_{,ij}M_k^{-\frac{4}{n-2}}
y^iy^j+O(|M_k^{-\frac 6{n-2}}|y|^3)\}v_k^2 \nonumber \\
&\ge & -CM_k^{2-\frac 8{n-2}}\int_{|y|\le \sigma M_k^{\frac{2}{n-2}}}
|y|^2|v_k^2-U^2|-
CM_k^{2-\frac {10}{n-2}}\int_{|y|\le \sigma M_k^{\frac{2}{n-2}}}
|y|^3U^2 \nonumber \\
&\ge &-C\sigma,\quad 3\le n\le 6.
\label{C2-1}
\end{eqnarray}

By (\ref{SS5}) and $|W|\to 0$ for $n=6$, 
\begin{eqnarray}
&&\limsup_{k\to \infty}M_k^2|I_3[v_k]| \nonumber \\
&=&\limsup_{k\to \infty} \frac{c(n)}2M_k^{2-\frac{4}{n-2}}
(\sigma M_k^{\frac 2{n-2}})\int_{|y|=\sigma M_k^{\frac 2{n-2}}}
\{\frac 12 R_{,ij}M_k^{-\frac 4{n-2}}y^iy^j+O(M_k^{-\frac 6{n-2}}|y|^3)\}v_k^2
\nonumber \\
&=&\limsup_{k\to \infty}\{\frac{c(n)}4\sigma M_k^{\frac{2(n-5)}{n-2}}
\int_{|y|=\sigma M_k^{\frac 2{n-2}}}R_{,ij}M_k^{-\frac 4{n-2}}y^iy^jU^2
\nonumber \\
&&+C\sigma M_k^{\frac{2(n-5)}{n-2}}
\int_{|y|=\sigma M_k^{\frac 2{n-2}}}|y|^2|v_k-U|U+
CM_k^{\frac{2(n-6)}{n-2}}\int_{|y|=\sigma M_k^{\frac 2{n-2}}}|y|^3U^2\}
\nonumber \\
&\le &C\sigma,\quad 3\le n\le 6.
\label{C3-1a}
\end{eqnarray}

Estimate (\ref{39-1new}) for $3\le n\le 6$ follows from 
(\ref{C1-1}), (\ref{C1-2}), (\ref{C2-1}) and (\ref{C3-1a}). To prove
(\ref{39-1new}) for $7\le n\le 9$, we make use of (\ref{2e-1}) with $\bar l=5$.

For $7\le n\le 9$, let 
$$\beta_2''=\beta_3'''=\left\{\begin{array}{ll}0,&\quad 7\le n\le 8,\\
1,&\quad n=9,
\end{array}
\right.
$$
$$\beta_2'''=n-6,\quad 7\le n\le 9.$$
By (\ref{nn1}), (\ref{nn2}) and (\ref{bb1}), (\ref{TT2}) 
holds for the above defined
$\beta_2''$, $\beta_3'''$ and $\beta_2'''$.
By (\ref{nn1}), (\ref{SS5}) and (\ref{second}),
\begin{eqnarray}
&&-M_k^2\frac{c(n)}2|S^{n-1}|\sum_{l=2}^5(l+2)M_k^{-\frac{4+2l}{n-2}}
\bar R^{(l)}\int_0^{\sigma M_k^{\frac 2{n-2}}}r^{l+n-1}U(r)^2dr \nonumber\\
&\ge &M_k^{\frac{2(n-8)}{n-2}}O(|W|)\int_0^{\sigma M_k^{\frac 2{n-2}}}
r^{3+n}U^2(r)dr \nonumber \\
&\to &0,\quad \mbox{as}\quad k\to \infty,\quad \mbox{for}\quad 7\le n\le 9.
\label{C5-1}
\end{eqnarray}
Multiplying (\ref{2e-1}) by $M_k^2$ with $\gamma_1=\frac 8{n-2}$,
$\gamma_2=n-6$, $\bar l=5$ and $\beta_2''$,$\beta_3'''$ given above, we obtain,
using also (\ref{C5-1}) that 
\begin{equation}
\liminf_{k\to \infty}M_k^2(I_5[v_k]+I_2[v_k])\ge 0.
\label{C6-1}
\end{equation}
We know 
\begin{eqnarray}
\lim_{k\to \infty}M_k^2I_3[v_k]&=&\lim_{k\to \infty}
\frac{\sigma}2c(n)\int_{|x|=\sigma}R(x)(M_ku_k)^2\nonumber \\
&=&\frac{\sigma}2c(n)\int_{|x|=\sigma}R(x)G^2.
\label{C6-2}
\end{eqnarray}
For the last term we are in conformal normal coordinates centered at $\bar P_1$.
It is elementary to see that for some constant $C_4(n)>0$,
\begin{equation}
G=C_4(n)r^{2-n}+O(r^{6-n})
\label{eq175}
\end{equation}

For $n=7$, $W(\bar P_1)=0$, so 
\begin{eqnarray}
&&\sigma \int_{|x|=\sigma}R(x)G^2 \nonumber \\
&=&\sigma \int_{|x|=\sigma}[\sum_{l=2}^3\sum_{|\alpha |=l}
\frac{\partial_{\alpha}R(\bar P_1)}{\alpha !}x^{\alpha}+O(|x|^4)]
\cdot [C_4(n)^2r^{4-2n}+O(r^{8-2n})]\nonumber \\
&=&\sigma \int_{|x|=\sigma}O(r^{8-2n})=O(\sigma),\quad n=7.
\label{C7-1}
\end{eqnarray}

For $n=8,9$, 
$$|R_{abcd}(\bar P_1)|=|R_{abcd,e}(\bar P_1)|=0.$$ Thus by
(\ref{second})
$$\int_{|x|=\sigma}\sum_{|\alpha |=4}\
\frac{\partial_{\alpha}R(\bar P_1)}{\alpha !}x^{\alpha}
=|S^{n-1}|\bar R^{(4)}(\bar P_1)=0.$$

It follows that
\begin{eqnarray}
&&\sigma \int_{|x|=\sigma}R(x)G^2 \nonumber \\
&=&\sigma \int_{|x|=\sigma}[\sum_{l=2}^5\sum_{|\alpha |=l}
\frac{\partial_{\alpha}R(\bar P_1)}{\alpha !}x^{\alpha}+O(|x|^6)]\cdot
[C_4(n)^2r^{4-2n}+O(r^{8-2n})]\nonumber \\
&=&\sigma \int_{|x|=\sigma}O(r^{10-2n})=O(\sigma),\quad n=8,9.
\label{C7-2}
\end{eqnarray}

By (\ref{C6-2}), (\ref{C7-1}), (\ref{C7-2}) and (\ref{C6-1})
$$\lim_{\sigma\to 0}\liminf_{k\to \infty}M_k^2I_5[v_k]\ge 0,\quad
7\le n\le 9. $$
Lemma \ref{Blem} is established. $\Box$

\bigskip

For $n=3,4,5$, 
$$G(\cdot,\bar P_1)=a_1(r^{2-n}+A+\mbox{higher order})$$
where $a_1>0$ and $A$ are constants. Thus, for some $\bar A\ge A$, 
$$G=a_1(r^{2-n}+\bar A+\mbox{higher order}).$$
For $n=6,7$, since $|W(\bar P_1)|=0$, we have
$$G(\cdot ,\bar P_1)=\left\{\begin{array}{ll}
a_2(r^{-4}+\psi(\theta)+O(r\log r),&\quad n=6,\\
a_3r^{-5}(1-a_4R_{,ij}(\bar P_1)x^ix^jr^2]+A+O(r),&\quad n=7
\end{array}
\right.
$$
where $x=r\theta$, $a_2,a_3>0$ and $A$ are constants. $\psi$ is a smooth
function on $\theta$. Thus
$$G=\left\{\begin{array}{ll}
a_2(r^{-4}+\bar \psi(\theta)+O(r\log r),&\quad n=6,\\
a_3r^{-5}(1-a_6R_{,ij}(\bar P_1)x^ix^jr^2)+\bar A+O(r),&\quad n=7,
\end{array}
\right.
$$
where $\bar A\ge A$ and $\bar \psi\ge \psi$.

A computation yields
\begin{equation}
\lim_{\sigma\to 0}B(\sigma, G, \nabla G)=
\left\{\begin{array}{ll}
-a_5\bar A,&\quad n=3,4,5,7\\
-a_6\int_{S^5}\bar \psi(\theta),&\quad n=6,
\end{array}
\right.
\label{C11-1}
\end{equation}
where $a_5,a_6>0$ are constants.

Consider on $M$,
$$\hat g=G(\cdot, \bar P_1)^{\frac 4{n-2}}g.$$
We make a change of variables $z=|x|^{-2}x$ and write 
$\hat g=\hat g_{ij}(z)dz^idz^j$. Since $W(\bar P_1)=0$ for $n=6,7$, it
is not difficult to verify that the mass of $(M,\hat g)$ is given,
modulo a positive constant multiple, by, see  \cite{bartnik},
$$\mbox{mass}=\lim_{\rho\to \infty}\int_{|z|=\rho}
(\partial_i\hat g_{ij}-\partial_j\hat g_{ii})\frac{z^i}{|z|}dz.$$
  
 A computation yields, modulo some positive constant multiple,
$$\mbox{mass}=\left\{\begin{array}{ll}
A,&\quad n=3,4,5,7\\
\int_{S^5}\psi,&\quad n=6.
\end{array}
\right.
$$
By the Positive Mass theorem of Schoen and Yau,
$A>0$ for $n=3,4,5,7$ and $\int_{S^5}\psi>0$ for $n=6$. Since 
$\bar A\ge A$ and $\bar \psi\ge \psi$, the above violates, in view of
(\ref{C11-1}) Lemma \ref{Blem}. Theorem \ref{thm3} is established. $\Box$ 

\bigskip

\noindent{\bf Comments on Remark \ref{rem0-1}:}  
For $n=8,9$, we not only have $|W(\bar x)|=0$, but also
$|\nabla W(\bar x)|=0$.  Since we work in conformal
local coordinates, $|\nabla R(\bar x)|=0$.
Then $G$ can be written as follows:  
$$G=\left\{\begin{array}{ll}
a_7(r^{-6}+\sum_{l=2}^4\psi_l(\theta)r^{l-4}+A+O(r\log r)),&\quad n=8,\\ \\
a_8(r^{-7}+\sum_{l=2}^4\psi_l(\theta)r^{l-5}+A+O(r)),&\quad n=9.
\end{array}
\right.
$$
where $a_7, a_8>0$ and
$\{\psi_l\}$ are spherical harmonics that have integral $0$ on $S^{n-1}$.
A contradiction can be obtained similarly from the positive mass theorem
for high dimensions.  Remark \ref{rem0-1} is proved. $\Box$

\section{A discussion of the case $1<1+\epsilon\le p\le \frac{n+2}{n-2}$}
The proofs for the  case $1<1+\epsilon\le p\le \frac{n+2}{n-2}$
are modifications of our proofs for $p=\frac{n+2}{n-2}$.
Modifications of similar nature can be found in \cite{SZ}, \cite{Li1} and
\cite{LiZhu}.
We 
point out some of these modifications for the proof of
Theorem \ref{thm1}.  We still prove it by
contradiction argument, so associated with $\{u_k\}$, we now have
$\{p_k\}\subset [1+\epsilon, \frac{n+2}{n-2}$.
Using standard blow up arguments together
with  the well known result in \cite{GS},
we only need to consider the case
 $\tau_k:=\frac{n+2}{n-2}-p_k\to 0$.
To avoid introducing new definitions we still let $v_k$ be defined as in 
(\ref{vk}). Then the equation
for $v_k$ is 
\begin{equation}
(\Delta +\bar b_i\partial_i+\bar d_{ij}\partial_{ij}-\bar c)v_k+n(n-2)
M_k^{-\tau_k}v_k^{p_k}=0,\quad |y|\le \delta M_k^{\frac{2}{n-2}}.
\label{vkp}
\end{equation}
where $\bar b_i$, $\bar d_{ij}$, $\bar c$ are defined as before.
Let $v_k^{\lambda}$ and $w_{\lambda}$ be defined as before, then direct
computation shows that
\begin{equation}
(\Delta +\bar b_i\partial_i+\bar d_{ij}\partial_{ij}-\bar c)w_{\lambda}
+n(n-2)p_kM_k^{-\tau_k}\xi_k^{p_k-1}w_{\lambda}\le E_{\lambda}
\quad \mbox{in}\quad \Sigma_{\lambda},
\label{etermp}
\end{equation}
where $E_{\lambda}$ and $\Sigma_{\lambda}$ are defined as before.
In this more general context, the estimate of 
$\sigma_k$ is related to that of $\tau_k$. First by Proposition \ref{prop2new}
and assuming $\sigma_k=\circ(1)$, we obtain
$$
E_\lambda=
 \bar c(y) U^\lambda(y)-
(\frac \lambda {|y|})^{n+2}\bar c(y^\lambda)U(y^\lambda)
+\circ(1)M_k^{ -\frac 4{n-2} } |y|^{-n}.
$$
Then we have 
$$\min_{|y|=L_k}v_k(y)\le (1+\epsilon)U(L_k),\quad L_k=
\left\{\begin{array}{ll}\delta M_k^{\frac{2}{n-2}},&\quad n=3,4,\\
M_k^{\frac{4}{(n-2)^2}},&\quad n\ge 5.
\end{array}
\right.
$$
where $\epsilon>0$ is an arbitrarily small positive constant, $\delta>0$ depends on 
$\epsilon$.
The above inequality leads to 
\begin{equation}
v_k(y)\le CU(y),\quad |y|\le \delta_1 L_k
\label{shpk}
\end{equation}
for some $\delta_1>0$ by the Moser iteration
technique used previously. The following Pohozaev identity will be used to 
obtain the vanishing rate of $\tau_k$:
\begin{eqnarray*}
&&\int_{|x|\le \sigma}\{(-b_i\partial_iu-d_{ij}\partial_{ij}u)(\nabla u\cdot x
+\frac{n-2}2u)
-\frac{c(n)}2u^2(x\cdot \nabla R)
-c(n)Ru^2\}\\
&&+\frac{\sigma}2c(n)\int_{|x|=\sigma}Ru^2
-\frac{n(n-2)}{p+1}\sigma\int_{|x|=\sigma}u^{p+1}
+(\frac{n^2(n-2)}{p+1}-\frac{n(n-2)^2}{2})\int_{B_{\sigma}}u^{p+1}\\
&=&B(\sigma,u,\nabla u)
\end{eqnarray*}
where
$B(\sigma,u,\nabla u)$ is defined as before.
Then by using (\ref{rough}) to evaluate the Pohozaev 
identity over $|y|\le L_k$ we obtain 
$\tau_k=O(M_k^{-2})$ for $n=3,4$ and 
$\tau_k=O(M_k^{-\frac{4}{n-2}})$ for $n\ge 5$. $\epsilon$ is arbitrarily small. Note that in 
the evaluation of Pohozaev identity, say $n\ge 5$, we first have 
$\tau_k=O(M_k^{-\frac{4}{n-2}+\tau_k})$, this implies
$M_k^{\tau_k}=O(1)$, consequently, $\tau_k=O(M_k^{-\frac{4}{n-2}})$.

For $n\ge 5$, by comparing the equations for $v_k$ and $U$ we know 
\begin{eqnarray}
&&(\Delta +\bar b_i\partial_i+\bar d_{ij}\partial_{ij}-\bar c)w_k
+n(n-2)M_k^{-\tau_k}p_k\xi_k^{p_k-1}w_k \nonumber \\
&=&\bar cU+n(n-2)(U^{\frac{n+2}{n-2}}-M_k^{-\tau_k}U^{p_k}),
\label{wkpk}
\end{eqnarray}
where $w_k=v_k-U$.
It is clear that for all $\epsilon>0$,
$$\bar cU+n(n-2)(U^{\frac{n+2}{n-2}}-M_k^{-\tau_k}U^{p_k})=
O(M_k^{-\frac{8}{n-2}})(1+r)^{4-n}+O(M_k^{-\frac{4}{n-2}+\epsilon})
(1+r)^{-2-n}.$$

Then by the Proof of Proposition \ref{prop5} we obtain 
$$\sigma_k\le CM_k^{-\frac{4}{n-2}+\epsilon},\quad \forall \epsilon>0.$$ 
Once we have this new estimate of $\sigma_k$ we can improve the estimate
of $E_{\lambda}$ to 
$$E_\lambda= O(1)M_k^{-\frac{8}{n-2}}|y|^{4-n}+O(1)M_k^{-\frac{8-\epsilon}{n-2}}|y|^{-n}.$$
So we can prove (\ref{shpk}) for a new $L_k$, which is
$$L_k:=\left\{\begin{array}{ll} \delta M_k^{\frac{2}{n-2}}&\quad n=5,\\
M_k^{\frac{8-\epsilon}{(n-2)^2}}&\quad n\ge 6.
\end{array}
\right.
$$
where $\epsilon>0$ is arbitrarily small. 

Now we want to use the new estimate of $\sigma_k$ to get a new estimate of
$\tau_k$. To do this, first for $n=5$ we compare $w_k$ with 
$$f(y)=QM_k^{-2+\frac{2\epsilon}{3}}|y|^{-\epsilon}
+QM_k^{-\frac{4-\epsilon}{3}}|y|^{\epsilon-3},\quad R<|y|<\delta M_k^{\frac{2}{3}}$$ 
where $\epsilon>0$ is small and $Q(\epsilon)>1$ is large, $R>>1$ is to make the 
maximum principle possible. Then from the maximum principle we get
$$|v_k(y)-U(y)|\le f(y),\quad \mbox{for}\quad R\le |y|\le \delta M_k^{\frac{2}{3}}.$$
Estimates for $|\nabla^j(v_k-U)|,j=1,2$ can be obtained similarly. Consequently we get
$$ v_k=U+O(M_k^{-2+\frac{2\epsilon}{3}})(1+|y|)^{-\epsilon}+O(M_k^{-\frac{4-\epsilon}{3}})
(1+|y|)^{\epsilon-3},\quad |y|<\delta M_k^{\frac{2}{n-2}}.$$
 Then we apply 
the Pohozaev identity over $|y|\le \delta M_k^{\frac 23}$ to get $\tau_k=O(M_k^{-2})$.

For $n\ge 6$, we follow the same procedure except that we use
$$f(y)=QM_k^{-\frac{8-2\epsilon}{n-2}}|y|^{-\frac{\epsilon}{10}}+
QM_k^{-\frac{4-\epsilon}{n-2}}|y|^{2-n+\epsilon},\quad R\le |y|\le 
M_k^{\frac{8-\sqrt{\epsilon}}{(n-2)^2}}.$$
Then the expansion for $v_k$ becomes
$$ v_k=U+O(M_k^{-\frac{8-2\epsilon}{n-2}})(1+|y|)^{-\frac{\epsilon}{10}}
+O(M_k^{-\frac{4-\epsilon}{n-2}})(1+|y|)^{2-n+\epsilon},
\quad |y|<M_k^{\frac{8-\sqrt{\epsilon}}{(n-2)^2}}.$$
 From the Pohozaev identity over $|y|\le M_k^{\frac{8-\sqrt{\epsilon}}{(n-2)^2}}$ 
we have $\tau_k=O(M_k^{-\frac{8-2\sqrt{\epsilon}}{n-2}})$, consequently
we obtain, by (\ref{wkpk}), that $\sigma_k=O(M_k^{-\frac{8-3\sqrt{\epsilon}}{n-2}})$.
>From now on we just replace $3\sqrt{\epsilon}$ by $\epsilon$ for convenience.
Now the right hand side of (\ref{wkpk}) becomes
$$\bar cU+n(n-2)(U^{\frac{n+2}{n-2}}-M_k^{-\tau_k}U^{p_k})=
O(M_k^{-\frac{8}{n-2}})(1+r)^{4-n}+O(M_k^{-\frac{8-\epsilon}{n-2}})
(1+r)^{-2-n}$$
Also the new estimate of $\sigma_k$ for $n\ge 6$ leads to a new estimate of 
$E_{\lambda}:$
$$ E_\lambda=\bar c(y) U^\lambda(y)-
(\frac \lambda {|y|})^{n+2}\bar c(y^\lambda)U(y^\lambda)
+O(1)M_k^{-\frac {12-\epsilon}{n-2}}|y|^{-n}.
$$ Then we obtain (\ref{shpk}) again
except that 
$$L_k=\left\{\begin{array}{ll}
\delta M_k^{\frac 2{n-2}},&\quad n=6,7 \\
M_k^{\frac {12-\epsilon}{(n-2)^2}},&\quad  n\ge 8.
\end{array}
\right.
$$

For $n=7$, we improve the estimate of $v_k-U$ by using 
$$f(y)=QM_k^{-\frac{8}{n-2}}|y|^{6-n}+QM_k^{-\frac{8-\epsilon}{n-2}}
|y|^{2-n+\epsilon},\quad R\le |y|\le \delta M_k^{\frac{2}{n-2}}$$
Then the Pohozaev identity gives $\tau_k=O(M_k^{-2})$. Note that in the 
computation, the term that contributes Weyl tensor and the term that 
contributes $\tau_k$ are both positive, therefore the existence of $W$ does
not affect the rate of $\tau_k$. After this we improve the rate of 
$\sigma_k$ to $O(M_k^{-\frac{8}{n-2}})$. 

For $n\ge 8$, to control $v_k-U$ over $L_k=M_k^{\frac {12-\epsilon}{(n-2)^2}}$
we use function 
$$f(y)=QM_k^{-\frac{8}{n-2}}|y|^{6-n+\bar a}+QM_k^{-\frac{8-\epsilon}{n-2}}
|y|^{2-n+\epsilon}$$
where $\bar a=\frac 23(n-8)+O(1)\epsilon$. Then from Pohozaev identity we get
$\tau_k=O(M_k^{-\frac{12-\epsilon}{n-2}})$. Then 
$\sigma_k=O(M_k^{-\frac{8}{n-2}})$ for $n\ge 8$. Once we have this rate of 
$\tau_k$, it does not affect the estimate of $\sigma_k$ any more (for $n=8,9$).
Then in the evaluation of the Pohozaev identity, the term that contributes $\tau_k$
has the same sign as the term that contributes all the curvature tensors. So
$\tau_k$ does not affect the vanishing rate of the Weyl or $\nabla R_{abcd}$ for 
$n=8,9$. We get the estimate of $\tau_k$ when we get the vanishing rate of 
curvature tensors.
Eventually we have $\tau_k=O(M_k^{-2})$ for $n=8,9$ and 
$\tau_k=O(M_k^{-\frac{16-\epsilon}{n-2}})$ for $n\ge 10$. Other estimates on 
$v_k$ (for $n=8,9$) are just the same as those in the special case 
$p_k\equiv \frac{n+2}{n-2}$. The proofs of Theorem \ref{thm3} and Theorem \ref{thm0}
can proceed as the special case $p_k\equiv \frac{n+2}{n-2}$.

\section{Appendix A: A useful function and its estimate}

For $n=3,4,5, \cdots$, let $A>2$, $a, b, \mu>0$,
and $0\le \gamma\le n-2$ be constants,
 we consider
functions $V$ and $H$ satisfying
\begin{equation}
-\mu^{-1}(1+r)^{-2-\mu} \le V(r)\le n(n+2)U(r)^{ \frac 4{n-2} }
+ \frac a2 r^{-4}, \qquad 1\le r\le A,
\label{V1}
\end{equation}
\begin{equation}
|V'(r)|\le \mu^{-1} r^{-3}, \qquad 1\le r\le A,
\label{NEW1}
\end{equation}
\begin{equation}
0\le H(r)\le b r^{\gamma-n},  \qquad 1\le r\le A,
\label{H1}
\end{equation}
and
\begin{equation}
|H'(r)|\le b r^{\gamma-n-1},\qquad 1\le r\le A,
\label{NEW2}
\end{equation}
where $U(r)=(1+r^2)^{ \frac {2-n}2 }$.
\begin{prop}  For $n=3,4,5, \cdots$,
let  $A>2$, $a, b>0$ and $0\le \gamma\le n-2$ be constants and
 let
$V$ satisfy (\ref{V1})-(\ref{NEW1}) and
$H$ satisfy (\ref{H1})-(\ref{NEW2}). Then there exists a unique
solution of
\begin{equation}
\left\{\begin{array}{ll}
\eta''(r)+\frac{n-1}r\eta'(r)
+(V(r)-\frac{a}{r^2})\eta(r)
=-H(r),\quad 1<r<A,\\ \\
\eta(1)=\eta(A)=0.
\end{array}
\right.
\label{eq1}
\end{equation}
Moreover
\begin{equation}
0\le \eta(r)\le Cr^{\gamma+2-n}, \qquad
1<r<A,
\label{eq2}
\end{equation}
and
\begin{equation}
|\eta'(r)|\le Cr^{ \gamma+1-n}, \quad
|\eta''(r)|\le Cr^{\gamma-n},\qquad
1\le r\le A,
\label{NEW3}
\end{equation}
where $C>0$ depends only
on $n, a$, $b$,  $\mu$ and $\gamma$.
\label{prop8}
\end{prop}

Thinking of $r=|x|$, $x\in \Bbb R^n$,  equation (\ref{eq1}) takes the
form
\begin{equation}
\left\{\begin{array}{ll}
\Delta \eta
+(V(|x|)-\frac{a}{|x|^2})\eta
=-H(|x|),\quad \mbox{in}\ B_A\setminus B_1,\\ \\
\eta=0,  \qquad\qquad
\qquad\qquad\qquad\qquad \mbox{on}\ \partial (B_A\setminus B_1).
\end{array}
\right.
\label{eq1new}
\end{equation}

\begin{lem}
For  $A>2$, let
 $\lambda_1=\lambda_1(n,A)$ be the first  eigenvalue of $-\Delta -
n(n+2)U^{ \frac 4{n-2} }$ on $B_A\setminus B_1$
with respect to zero Dirichlet boundary value.  Then
$\lambda_1>0$.
\label{lem7-0}
\end{lem}

\noindent{\bf Proof of Lemma \ref{lem7-0}.}\
Since $U_t(x):=t^{ \frac {n-2}2 }U(tx)
$ satisfies, for all $t>0$,
$$
-\Delta U_t-n(n-2)U_t^{ \frac {n+2}{n-2} }=0, \qquad
\mbox{in}\ \Bbb R^n,
$$
$\frac {d}{dt}U_t|_{ t=1}$ satisfies the linearized equation, i.e.,
$$
-\Delta \varphi - n(n+2) U^{ \frac 4{n-2} }\varphi=0,
\qquad
\mbox{in}\ \Bbb R^n,
$$
where $\varphi(r):=\frac {r^2-1}{ r^2+1}U(r)$ is
positive in   $1<r<\infty$.

Let $\bar \eta$ be a positive eigenfunction with respect
to $\lambda_1$, so
$$
-\Delta \bar \eta-n(n+2)U^{ \frac 4{n-2} }
\bar\eta=\lambda_1\bar\eta,
\qquad 1<|x|<A,
$$
and
$
\bar\eta=0$ on $\partial (B_A\setminus B_1)$.
Multiplying the above equation of $\bar \eta$  by $\varphi$
and integrating by parts lead to
$$
\lambda_1\int_{ B_A\setminus B_1}
\bar\eta \varphi> \int_{ B_A\setminus B_1}
[- \bar\eta\Delta \varphi-
n(n+2)U^{ \frac 4{n-2} } \bar\eta \varphi]
=0.
$$
Lemma \ref{lem7-0} is  established.

\vskip 5pt
\hfill $\Box$
\vskip 5pt

\begin{cor} Under the hypotheses of
Proposition \ref{prop8},
 equation (\ref{eq1}) has a unique solution $\eta$, which is
non-negative.
\label{cor7-1}
\end{cor}

\noindent{\bf Proof of Corollary \ref{cor7-1}.}\
By (\ref{V1}),
\begin{equation}
V(|x|)-\frac{a}{|x|^2}\le n(n+2)U(x)^{ \frac 4{n-2} }
-\frac{a}{2|x|^2},
\qquad 1 < |x|<A.
\label{be2}
\end{equation}
Corollary \ref{cor7-1} follows from
 Lemma
\ref{lem7-0} and standard elliptic theories.

\noindent{\bf Proof of Proposition \ref{prop8}.}\
Fix an $R>2$, depending only on $n$ and $a$, such that
\begin{equation}
 n(n+2)U(r)^{ \frac 4{n-2} }
+\frac a2 r^{-4}-a r^{-2}
\le -\frac a2 r^{-2}, \qquad \forall\ r\ge R.
\label{RR12}
\end{equation}
If $A\le 3R$, we know from  Lemma \ref{lem7-0} that
 the first
eigenvalue $\lambda_1$ of $-\Delta -(V-\frac a{r^2})$
on $B_A\setminus B_1$, with respect to
zero Dirichlet boundary value,
has a positive lower bound
which depends only on
$n$ and $R$.
Thus  the $L^2$
norm of $\eta$ on $B_A\setminus B_1$ is under control.
By standard elliptic estimates, the $L^\infty$ norm of
$\eta$  on $B_A\setminus B_1$ is also under control.
In the following we assume that 
$A>3R$.
By (\ref{V1}) and (\ref{RR12}),
\begin{equation}
V(r)-ar^{-2}\le -\frac a2 r^{-2}, 
\qquad \forall R\le r\le A.
\label{RR13}
\end{equation}
Since $0<\gamma \le n-2$, we can pick
 some constant $C>1$, depending only on
$\gamma, n, b, a$ such that
$w(r):= C r^{ \gamma+2-n}$ satisfies 
\begin{equation}
[\Delta +(V-a r^{-2})]w\le
-b r^{\gamma -n}
\le -H(r), \qquad
\forall\ R\le r\le A.
\label{RR14}
\end{equation}

Fix some smooth function $f(r)$, depending only on
the usual parameters (i.e. $n, a, b, \mu$ and $\gamma$), satisfying
$$
f(r)\equiv 0, \qquad 2R<r<\infty,
$$
$$
f(r)\le -[\Delta +(V-a r^{-2})]w
-H(r), \qquad 1\le r\le R,
$$
and
$$
f(r)\le 0, \qquad \forall\ 1\le r<\infty.
$$

To prove (\ref{eq2}), we
only need to find some non-negative function $w_1(r)$ 
in $1\le r<\infty$ satisfying
\begin{equation}
w_1(r)\le C_1 r^{2+\gamma-n}, \qquad
1\le r<\infty,
\label{21-1new}
\end{equation}
and
\begin{equation}
[\Delta +(V-ar^{-2})]w_1(r)\le f(r),
\qquad \forall\ 1\le r\le A,
\label{21-2new}
\end{equation}
where $C_1>0$ is some constant depending only on the
usual parameters.

Indeed, let $w_1$ be as above, then
$$
[\Delta +(V-ar^{-2})](w+w_1)\le -H_1(r), \qquad
1\le r\le A,
$$
and therefore, in view of (\ref{eq1}),
$$
[\Delta +(V-ar^{-2})](w+w_1-\eta)\le 0,
\qquad \mbox{on}\ B_A\setminus B_1.
$$
We also have, using the non-negativity of $w$ and $w_1$,
$$
w+w_1-\eta\ge 0, \qquad \mbox{on}\ \partial (B_A\setminus B_1).
$$
Because of Lemma \ref{lem7-0} and (\ref{be2}), we may apply the
maximum principle to obtain
$$
 w+w_1-\eta\ge 0, \qquad 1\le r\le A.
$$
This  gives the desired  estimate  (\ref{eq2}).

Now  we construct  such a $w_1$. 
Consider 
\begin{equation}
\left\{
\begin{array}{rll}
\Delta \tilde w_1(y)
+[n(n+2) U(y) ^{ \frac 4{n-2} }-\frac a2]\tilde w_1(y)
&=& |y|^{ -n-2} f(\frac y {|y|^2}), 
\qquad |y|<1,\\
\tilde w_1(y)&=&0,\qquad |y|=1.
\end{array}
\right.
\label{ww1}
\end{equation}

We know from the proof of Lemma \ref{lem7-0} that
$\varphi(r):=\frac {r^2-1}{r^2+1} U(r)$
satisfies
$$
[\Delta +n(n+2)U(y)^{ \frac 4{n-2} }]\varphi(y)=0,
\qquad |y|<1,
$$
$$
\varphi(r)<0, \quad \forall\ 0\le r<1, \
\qquad \mbox{and}\qquad \varphi(1)=0.
$$
It follows that the first eigenvalue, with respect
to zero Dirichlet boundary data,  of 
$-\Delta - n(n+2) U^{ \frac 4{n-2} }$ on $B_1$ is zero.
So the first eigenvalue of
$-\Delta - n(n+2) U^{ \frac 4{n-2} } +\frac a2$ on $B_1$ is 
equal to $\frac a2>0$.
We also know that $|y|^{-n-2}f(\frac y{|y|^2})$ is
non-positive for all $|y|\le 1$ and is equal to
zero for $|y|\le \frac 1{2R}$.
By standard elliptic theories, (\ref{ww1}) has a unique
radial  solution
$\tilde w_1$ satisfying
$$
0\le \tilde w_1(y)\le C, \qquad \forall\ |y|\le 1, 
$$
where  $C$ is some positive constant  depending only on the usual parameters.

Let 
$$
w_1(x)=\frac 1 {|x|^{n-2} } \tilde w_1( \frac x { |x|^2 } ),
\qquad |x|\ge 1.
$$
Then, because of  (\ref{ww1}), $w_1$ satisfies
$$
\big\{\Delta +[  n(n+2) U(r)^{ \frac 4{n-2} }
-\frac a{ 2 r^4}]\big\} w_1(r)=f,   
\qquad
r> 1,
$$
and 
$$
0\le w_1(r)\le Cr^{2-n}, \qquad \forall r\ge 1.
$$
Finally, using  (\ref{be2}),
we have, for $1\le r\le A$,
\begin{eqnarray*}
[\Delta+(V(r)-ar^{-2})]w_1
&\le &
[\Delta +(n(n+2)U(r)^{4/(n-2)}- \frac a{2r^2})]w_1\\
&\le &[\Delta +(n(n+2)U(r)^{4/(n-2)}- \frac a {2r^4})]w_1
=f.
\end{eqnarray*}
Thus, the $w_1$ has the desired properties, 
and (\ref{eq2}) is established.
Using (\ref{NEW1}), (\ref{NEW2}), (\ref{eq2}) and 
(\ref{eq1}), estimate
(\ref{NEW3}) follows from standard elliptic theories
with the help of a standard scaling argument.
Proposition \ref{prop8} is established.

\vskip 5pt
\hfill $\Box$

yyli@math.rutgers.edu  \qquad leizhang@math.ufl.edu
\end{document}